\documentclass[twoside,a4paper,12pt]{article}
\usepackage{amsmath}
\usepackage[usenames,dvipsnames]{color}

\def\juerg #1{{\color{red}#1}}
\let\juerg\relax

\newcommand{\dega}{{\Delta_\Gamma}}
\newcommand{\buuga}{{(\bar u, \bar u_\Gamma)}}
\newcommand{\nf}{{\bf n}}
\newcommand{\uga}{{u_\Gamma}}
\newcommand{\uuga}{{(u,u_\Gamma)}}
\newcommand{\yga}{{y_\Gamma}}
\newcommand{\vga}{{v_\Gamma}}
\newcommand{\rz}{{\rm I\!R}}
\newcommand{\nz}{{\rm I\!N}}

\newcommand{\dx}{{{\rm d}x}}
\newcommand{\dt}{{{\rm d}t}}
\newcommand{\ds}{{{\rm d}s}}
\newcommand{\dgm}{{{\rm d}\Gamma}}
\newcommand{\essinf}{\mathop{\rm ess\, inf}}
\newcommand{\esssup}{\mathop{\rm ess\, sup}}
\newcommand{\oma}{{\Omega}}
\newcommand{{\tinto}}{{\int_0^T}}
\newcommand{{\xinto}}{{\int_\Omega}}
\newcommand{{\txinto}}{{\int_0^T\!\!\int_\Omega}}
\newcommand{{\tgamma}}{{\int_0^T\!\!\int_\Gamma}}
\newcommand{\txt}{{\int_0^t\!\!\int_\Omega}}
\newcommand{\txg}{{\int_0^t\!\!\int_\Gamma}}
\newcommand{\lzo}{{L^2(\Omega)}}
\newcommand{\heins}{{H^1(\Omega)}}
\newcommand{\lio}{{L^\infty(\Omega)}}

\newcommand{\gheins}{{H^1(\Gamma)}}
\newcommand{\qlzo}{{L^2(Q)}}
\newcommand{\qlio}{{L^\infty(Q)}}
\newcommand{\glzsig}{{L^2(\Sigma)}}
\newcommand{\glisig}{{L^\infty(\Sigma)}}
\newcommand{\uad}{{\cal U}_{\rm ad}}
\newcommand{\qed}{\hfill\colorbox{black}{\hspace{-0.01cm}}}

\title{{\Huge Optimal control of an Allen--Cahn\\[2mm]
equation with singular potentials\\[2mm]
and dynamic boundary condition\\[2mm]}}
\date{}

\author{Pierluigi Colli\thanks{Dipartimento di Matematica  ``F. Casorati'',
Universit\`a di Pavia, Via Ferrata, 1, 27100 Pavia,  Italy 
({\tt pierluigi.colli@unipv.it}).}
        \and J\"urgen Sprekels\thanks{Weierstrass Institute for 
Applied Analysis and Stochastics,
Mohrenstrasse 39, 10117 Berlin, Germany
 ({\tt juergen.sprekels@wias-berlin.de}).}}

\begin{document}

\maketitle

\begin{abstract}
In this paper, we investigate optimal control problems for Allen--Cahn
equations with 
\juerg{differentiable}
singular nonlinearities and a dynamic boundary
condition involving 
\juerg{differentiable}
singular nonlinearities and the Laplace--Beltrami
operator. The approach covers both the cases of distributed controls
and of boundary controls. The cost functional is of standard tracking
type, and box constraints for the controls are prescribed.  Parabolic
problems with nonlinear dynamic boundary conditions involving the
Laplace--Beltrami operator have recently drawn increasing attention
due to their importance in applications, while their optimal control
was apparently never studied before. In this paper, we first extend
known well-posedness and regularity results for the state equation and
then show the existence of optimal controls and that the
control-to-state mapping is twice continuously Fr\'echet
differentiable between appropriate function spaces. Based on these
results, we establish the first-order necessary optimality conditions
in terms of a variational inequality and the adjoint state equation,
and we prove second-order sufficient optimality conditions.\\[2mm]
\end{abstract}

\pagestyle{myheadings}
\thispagestyle{plain}
\markboth{\sc Optimal control of an Allen--Cahn equation}{\sc
  Pierluigi Colli and J\"urgen Sprekels}

\section{Introduction}

Let $\oma\subset\rz^N$, $2\le N\le 3$, denote some open and bounded
domain with smooth boundary $\Gamma$ and outward unit normal $\nf$, and let $T>0$ be a
given final time. We put $Q:=\oma\times (0,T)$ and $\Sigma:=\Gamma\times (0,T)$,
and we assume that $\beta_i\ge 0$, $1\le i\le 6$, are given
constants which do not all vanish. Moreover, we assume:

\vspace{3mm}
{\bf (A1)}\quad There are given functions
\begin{eqnarray*}
&&z_Q\in L^2(Q),
\quad z_\Sigma\in L^2(\Sigma),\quad z_T\in H^1(\oma),
\quad z_{\Gamma,T}\in H^1(\Gamma),\\[1mm]
&&\widetilde{u}_1, \widetilde{u}_2\in L^\infty(Q) \mbox{\,\, with \,\,}\widetilde{u}_1\le \widetilde{u}_2
\mbox{\,\, a.\,e. in \,}\,Q,\nonumber\\[1mm]
&&\widetilde{u}_{1_\Gamma}, \widetilde{u}_{2_\Gamma}\in L^\infty(\Sigma)
\mbox{\,\, with \,}\, \widetilde{u}_{1_\Gamma}\le 
\widetilde{u}_{2_\Gamma} \mbox{\,\, a.\,e. on \,\,} \Sigma\,.
\end{eqnarray*}

 We then consider the following (tracking type) optimal control problem: 

{\bf (CP)}\quad Minimize
\begin{eqnarray}
\label{eq:1.1}
 J((y,\yga),(u,\uga))&:=&\frac {\beta_1} 2 \tinto\!\!\xinto \left 
|y-z_Q\right|^2\,\dx\,\dt\,+\,\frac {\beta_2}
2\tinto\!\!\int_\Gamma\left|\yga-z_\Sigma\right|^2\,\dgm \,\dt
\nonumber\\[1mm] 
&&+\,\frac{\beta_3}2 \xinto\left|y(\cdot,T)-z_T\right|^2\dx\, +
\,\frac {\beta_4} 2\int_\Gamma \left|\yga(\cdot,T)-z_{\Gamma,T}
\right|^2\, \dgm\nonumber\\[1mm]
&&+\,\frac{\beta_5}2 \tinto\!\!\xinto |u|^2\,\dx\,\dt\,+\,\frac {\beta_6} 2
\tinto\!\!\int_\Gamma \left|\uga\right|^2
\,\dgm\,\dt
\end{eqnarray}
subject to the parabolic initial-boundary value problem with nonlinear dynamic boundary condition
\begin{eqnarray}
\label{eq:1.2}
&y_t-\Delta y+f'(y)=u\,\quad\mbox{a.\,e. in }\,Q,&\\[2mm]
\label{eq:1.3}
&\partial_t y_{\Gamma}-\dega \yga + \partial_\nf y + g'(\yga)=\uga,
\quad y_{|\Gamma}=\yga, \,\quad\mbox{a.\,e. on }\,
\Sigma,&\\[2mm]
\label{eq:1.4}
&y(\cdot,0)=y_0 \,\quad\mbox{a.\,e. in }\,\oma,\,\quad\,y_\Gamma(\cdot,0)=y_{0_\Gamma}\,\quad
\mbox{a.\,e. on }\,\Gamma,&
\end{eqnarray}
and to the control constraints
\begin{eqnarray}
&&(u,\uga)\in \uad := \left\{(w,w_\Gamma)\in L^2(Q)\times L^2(\Sigma)\,:\,
\widetilde{u}_1\le w\le \widetilde{u}_2\quad 
{\hbox{a.\,e. in } Q},\right.\nonumber\\[1mm]
&&\hskip4cm  \left. \widetilde{u}_{1_\Gamma}\le w_\Gamma\le
\widetilde{u}_{2_\Gamma} \quad {\mbox{a.\,e. in }\Sigma} \ \right\}.
 \label{eq:1.5}
\end{eqnarray}
Here,  $y_0$ and $y_{0_\Gamma}$ are given initial data 
\juerg{such that $y_{0_\Gamma}=y_{0|\Gamma}$}, 
$\dega$ is the Laplace--Beltrami 
operator on $\Gamma$, and the functions $f\,,\,g$ are given nonlinearities,
while 
\juerg{$u$ and $\uga$} 
play the roles of distributed 
\juerg{and}
boundary controls, respectively. Note that we do not require $\uga$ to be {somehow} the 
restriction of $u$ on $\Gamma$; such a requirement would be much too restrictive 
for a control to satisfy. 

We remark at this place that for the cost functional to be meaningful it would suffice to only
assume that $z_T\in \lzo$ and $z_{\Gamma,T}\in L^2(\Gamma)$. However, the higher regularity
of $z_T$ and $z_{\Gamma_T}$ requested in {\bf (A1)} will later be essential to be able to
treat the adjoint state problem. 

The system (\ref{eq:1.2})--(\ref{eq:1.4}) is an initial-boundary 
value problem with nonlinear dynamic boundary condition for an Allen--Cahn equation.
In this connection, the unknown $y$ usually stands for the order parameter of an
isothermal phase transition, typically the fraction of one of the involved phases. 
In such a situation
it is physically meaningful to require $y$ to attain values in the interval $[0,1]$ on
both $\oma$ and $\Gamma$. A standard technique to meet this requirement is to postulate
that the first derivatives of the bulk potential $f$ and of the surface potential $g$
become singular at 0 and at 1. A typical form for such a potential is $f=f_1+f_2$, 
where $f_2$ is smooth on $[0,1]$ and $f_1(y)= \alpha\, [y\ln(y)+(1-y)\ln(1-y)]$ with some
$\alpha>0$. Another possibility is to choose $f_1$ as the indicator function $I_{[0,1]}$ of 
the interval $[0,1]$; in this case $f_1'$ has to be replaced by the subdifferential
$\partial I_{[0,1]}$, and $f$ becomes a {\em double obstacle} potential. In this case,
(\ref{eq:1.2}) has to be understood as a differential inclusion or variational inequality. 
Similar choices can be made for the surface potential $g$. 
\juerg{In this paper, we confine ourselves to differentiable singularities such as
the logarithmic type. The non-differentiable case of indicator functions, which requires quite different techniques, will be the
subject of the forthcoming paper \cite{CFS}; in this nonsmooth case only results concerning
existence and first-order necessary optimality conditions can be derived while
the question of sufficient optimality conditions remains unanswered.}

\juerg{For a discussion of the physical background of the dynamical boundary condition (\ref{eq:1.3}), we refer the reader to \cite{WS} and the literature cited therein. In particular, it is related to some phase transition problems in materials science. For instance, physicists have pointed out 
(cf. \cite{KEM} and the references given there) that in connection with phase separation phenomena for certain materials 
a dynamical interaction
with the wall (i.e., $\Gamma$) must be taken into account. This fact corresponds to considering a free energy functional that also contains a
boundary contribution. As a consequence, they deduced a (linear) dynamic boundary
condition of the form
\begin{equation}
\label{eq:1.3n}
\partial_t y_{\Gamma}-\dega \yga + \partial_\nf y + \yga= 0\,.
\end{equation}
Phenomenologically speaking, the boundary condition (\ref{eq:1.3n}) means that 
on the surface the order parameter 
 relaxes towards equilibrium at a rate proportional to the driving
force given by the Fr\'echet derivative of the free energy functional. The more
general condition (\ref{eq:1.3}) guarantees that $y_{|\Gamma}$ attains 
values in $[0,1]$, and the boundary control $\,\uga\,$ reflects the possibility
of monitoring the phase separation processes  by
supplying or removing mass at the boundary.}

There exists a vast literature on the well-posedness and asymptotic behaviour
of the Allen--Cahn equation with the {\em no-flux boundary condition}
$\partial_\nf y=0$ in place of (\ref{eq:1.3}). Also, the well-posedness and asymptotic
behavior of the system (\ref{eq:1.2})--(\ref{eq:1.4}) has been the subject of 
numerous papers (see \cite{CC} and the many references given there). 

Moreover, 
distributed and boundary control problems for the Allen--Cahn equation with no-flux 
boundary conditions or boundary conditions of the third kind have been studied in a number 
of recent papers, in particular, for the case of the double obstacle potential. In this
connection, we refer to \cite{FS1} and \cite{FS2}. Associated stationary, that is,
elliptic MPEC problems have been studied in \cite{HS} (see also the monograph \cite{NST}),
and the related Cahn--Hilliard case was recently analyzed in \cite{HW}. We also like
mention the works \cite{Hei, HeiTr, HJ} that treated optimal control problems 
for the Caginalp-type temperature-dependent generalization of the Allen--Cahn equation in the case of 
nonsingular potentials and standard boundary conditions; a thermodynamically consistent
temperature-dependent model with singular potential of the above logarithmic type 
was the subject of \cite{LS}.

The main novelty of the present paper is to study optimal control problems with singular 
potentials of the logarithmic type and dynamic boundary conditions of the form 
(\ref{eq:1.3}).  
In fact, while various types of dynamic boundary conditions have already been studied in connection with 
optimal control theory (see~\cite{HKR}, for a recent example), it seems that dynamic
boundary conditions involving the Laplace--Beltrami operator have not been considered before.
One of the difficulties is that from the viewpoint of optimal control it does not make
sense to postulate that the controls $u$ and $\uga$ satisfy $u_{|\Gamma}=\uga$.   

The paper is organized as follows: in Section 2, we give a precise statement of the problem
under investigation, and we derive some results concerning the state system (\ref{eq:1.2})--(\ref{eq:1.4})
and a certain linear counterpart, which will be employed repeatedly in the later analysis. In Section 3, we then treat the optimal control problem,
proving the existence of optimal controls and deriving the first-order necessary and the
second-order sufficient optimality conditions. During the course of this analysis, we will make 
repeated use of the elementary Young's inequality
$$
a\,b\,\le\,\gamma |a|^2\,+\,\frac 1{4\gamma}\,|b|^2\quad\forall\,a,b\in\rz \quad\forall\,\gamma>0,
$$
and of the fact that we have the continuous embeddings $\heins\subset L^p(\oma)$,
for $1\le p\le 6$, and $H^2(\oma)\subset\lio$ in three dimensions of space. In particular, we have
\begin{eqnarray}
\label{eq:1.6}
&&\|v\|_{L^p(\oma)}\,\,\le\,\widetilde C_p\,\,\|v\|_{\heins}\,\,\,
\quad\forall \,v\in \heins, \quad 1\le p\le 6,\\[1mm]
&&\|v\|_\lio\,\le\,\widetilde C_\infty\,\|v\|_{H^2(\oma)}\,\quad\forall\,v\in H^2(\oma),
\end{eqnarray}
with positive constants $\widetilde C_p$ that only depend on $\oma$.

\section{General assumptions and the state equation}
\setcounter{equation}{0}
In this section, we formulate the general assumptions of the paper, and we state some results for the state system (\ref{eq:1.2})-(\ref{eq:1.4}). To this end, we introduce the function spaces
\begin{eqnarray}
\label{eq:2.1}
&&H:=L^2(\oma),\quad V:=H^1(\oma),\quad H_\Gamma:=L^2(\Gamma), \quad V_\Gamma:=H^1(\Gamma),\nonumber\\[2mm]
&&{\cal H}:=\qlzo\times\glzsig,\quad {\cal X}:=\qlio\times\glisig,\nonumber\\[2mm]
&&{\cal Y}:=\left\{(y,y_\Gamma):\,y\in H^1(0,T;H)\cap C^0([0,T];V)\cap L^2(0,T;H^2(\Omega)),\right.\nonumber\\[1mm]
&&\quad\quad\quad \left.y_\Gamma \in H^1(0,T;H_\Gamma)\cap C^0([0,T];V_\Gamma)\cap L^2(0,T;H^2(\Gamma)),\quad
y_\Gamma=y_{|\Gamma}    \right\},\quad\quad
\end{eqnarray}

\noindent
which are Banach spaces when endowed with their natural norms. In the following, we denote the norm in a Banach
space $E$ by $\|\,\cdot\,\|_E$; for convenience, the norm of the space $H\times H\times H$ will also be denoted by
$\|\,\cdot\,\|_H$.
Identifying $H$ with its dual space $H^*$, we have the Hilbert triplet $V\subset H\subset V^*$, with dense and compact embeddings. Analogously, we obtain the triplet $V_\Gamma\subset H_\Gamma\subset V_\Gamma^*$, with dense and compact embeddings. 
We make the following general assumptions: 

\vspace{3mm}
{\bf (A2)}\quad$f=f_1+f_2$ and $g=g_1+g_2$, where $f_2, g_2\in C^3[0,1]$, \,and where \,$f_1,g_1\in C^3(0,1)$\,
are convex and satisfy the following conditions:
\begin{eqnarray}
\label{eq:2.3}
&&\lim_{r\searrow 0} \,f_1'(r)=\lim_{r\searrow 0} \,g_1'(r)=-\infty\,,\quad
\lim_{r\nearrow 1} \,f_1'(r)=\lim_{r\nearrow 1} \,g_1'(r)=+\infty\,.\\[2mm]
\label{eq:2.4}
&&\exists \,\,M_1\ge 0, \,M_2>0\,\quad\mbox{such that }\,\,\,|f_1'(r)|\le 
M_1+M_2\,|g_1'(r)|\,\quad\forall\,r\in (0,1).\quad
\end{eqnarray}

\vspace{3mm}
{\bf (A3)}\quad$y_0\in V$, \,\,$y_{0_\Gamma}\in V_\Gamma$, and we have 
$f_1(y_0)\in L^1(\oma)$, $g_1(y_{0_\Gamma})\in L^1(\Gamma)$, and
\begin{equation}
\label{eq:2.2}
0<y_0<1\quad\mbox{a.\,e. in }\,\oma,\,\,\quad\,\,
0<y_{0_\Gamma}<1 \quad\mbox{a.\,e. on }\,\Gamma\,.
\end{equation}
 
{\bf Remark 1:} \,\,The condition (\ref{eq:2.3}) is obviously satisfied if $f_1$ and $g_1$ are potentials 
of logarithmic type as those mentioned in Section 1, while (\ref{eq:2.4}) is needed for the existence result from \cite{CC} that will be used below.

\vspace{2mm} 
To simplify notation, in the following we will denote the trace $y_{|\Gamma}$ (if it exists) of a function $y$ on $\Gamma$ by
$y_\Gamma$ without further comment. Now observe that the set $\uad$ is a bounded subset of ${\cal X}$. Hence, there exists
a bounded open ball in ${\cal X}$ that contains $\uad$. For later use it is convenient to fix such a ball once
and for all, noting that any other such ball could be used instead. In this sense, the following assumption
is rather a denotation:

\vspace{2mm}
{\bf (A4)}\quad${\cal U}$ is a nonempty open and bounded subset of ${\cal X}$ containing $\uad$, and the constant
$R>0$ satisfies
\begin{equation}
\label{eq:2.5}
\|u\|_\qlio\,+\,\|u_\Gamma\|_\glisig\,\le\,R \quad\,\forall\,(u,u_\Gamma)\in {\cal U}.
\end{equation}  

\vspace{2mm}
The following result is a special case of \cite[{Theorems~2.3--2.5} and Remark~4.5]{CC} if one
puts (in the notation of \cite{CC}) $\,\beta=f_1'$, $\beta_\Gamma=g_1'$, $\pi=f_2'$, $\pi_\Gamma=
\juerg{g_2'}$ there.

\vspace{3mm}

{\bf Theorem~2.1} \,\,{\em Suppose that the general assumptions} {\bf (A2)}, {\bf (A3)} {\em are satisfied.
Then we have:}

(i)  \,\,\juerg{{\em For any pair} $(u,\uga)\in{\cal H}$, {\em the state system} (\ref{eq:1.2})--(\ref{eq:1.4}) 
{\em has a unique solution $(y,y_\Gamma)\in {\cal Y}$ such that}}
$ 0<y <1\quad\mbox{{\em a.\,e. in }}\,Q\,\mbox{{\em and }}\,
0<
\juerg{y_\Gamma}
<1 \quad\mbox{{\em a.\,e. on }}\,\Sigma\,.
$

(ii) \, {\em Suppose that also} {\bf (A4)} {\em is fulfilled. 
Then there is a positive constant $K^*_1$, which only
depends on $\oma$, $T$, $y_0$, $y_{0_\Gamma}$, $f$, $g$, and $R$, such that for every $(u,\uga)\in {\cal U}$ the associated solution $(y,y_\Gamma)\in {\cal Y}$ satisfies
\begin{eqnarray}
\|(y,y_\Gamma)\|_{\cal Y}&=&\|y\|_{H^1(0,T;H)\cap C^0([0,T];V)\cap L^2(0,T;H^2(\Omega))}
\nonumber\\[1mm]
&&+\,\|y_\Gamma\|_{H^1(0,T;H_\Gamma)\cap C^0([0,T];V_\Gamma)\cap L^2(0,T;H^2(\Gamma))}\,\,\le\,\,
K_1^*\,, \label{eq:2.6} \\[2mm]
&& \hskip-1cm { \|f'(y) \|_{L^2(0,T;H)}\, +\,\|g'(y_\Gamma)\|_{L^2(0,T;H_\Gamma)} \,\,\le\,\,
K_1^*\,.} \label{eq:2.6bis}
\end{eqnarray}
Moreover, there is some $K^*_2>0$, which only
depends on $\oma\,,\, T\,,\, y_0\,,\, y_{0_\Gamma}\,,\, f\,,\, g\,,
\mbox{{\em and }} R\,,$ such
that we have: whenever $\!(u_1, u_{1_\Gamma}), 
(u_2, u_{2_\Gamma}) \in {\cal U}\!$ are given and 
$(y_1,y_{1_\Gamma}), (y_2,y_{2_\Gamma})\in {\cal Y}$  
denote the associated solutions of the state system, then we have}
\begin{eqnarray}
\label{eq:2.7}
\hspace*{-8mm}&&\|y_1-y_2\|_{C^0([0,T];H)}^2\,+\,\|\nabla (y_1-y_2)\|_{L^2(Q)}^2 \,+\,
\|y_{1_\Gamma}-y_{2_\Gamma}\|_{C^0([0,T];H_\Gamma)}^2\,\nonumber\\
\hspace*{-8mm}&& +\,\|\nabla_\Gamma(y_{1_\Gamma}-y_{2_\Gamma})\|_{L^2(\Sigma)}^2
\le\,K_2^*\,\left\{\|u_1-u_2\|_{L^2(0,T;H)}^2\,+\,\|u_{1_\Gamma}-u_{2_\Gamma}\|^2_{L^2(0,T;H_\Gamma)}
\right\}\!.
\end{eqnarray}  

\vspace{2mm}
{\bf Remark 2:} \,\,(i)\, It follows from Theorem~2.1, in particular, that the control-to-state mapping   
${\cal S}$, $(u,\uga)\mapsto {\cal S}(u,\uga):=(y,\yga)$ is well defined as a mapping from 
${\cal X}$ into ${\cal Y}$; moreover, ${\cal S}$
is Lipschitz continuous when viewed as a mapping from the subset ${\cal U}$ of ${\cal H}$
into the space 
$$\bigl(C^0([0,T];H)\cap L^2(0,T;V)\bigr)\times
\bigl(C^0([0,T];H_\Gamma)\cap L^2(0,T;V_\Gamma)\bigr).$$

\vspace{2mm} (ii)\, Observe that we cannot expect $y$ to be continuous in $\overline Q$, since both
$\partial_t y_\Gamma$ and $\dega \yga$ only belong to $\glzsig$, so that also only $\partial_\nf\yga\in\glzsig$. 
However, we have $y\in L^2(0,T;C^0(\overline \oma))$ and $\yga\in L^2(0,T;C^0(\Gamma))$.

\vspace{3mm}
The next result is concerned with a linear problem with dynamic boundary condition. It will later be needed
to ensure the solvability of a number of linearized systems.

\vspace{2mm}
{\bf Theorem~2.2} \,\,{\em Suppose that functions $(u,u_\Gamma)\in {\cal H}$, $c_1\in\qlio$,
$c_2\in\glisig$, 
\juerg{and $w_0\in\heins$ with $w_{0_\Gamma}:=w_{0|\Gamma}\in H^1(\Gamma)$}
are given. Then we have:}

(i) \,\,{\em The linear initial-boundary value problem}
\begin{eqnarray}
\label{eq:2.8}
&w_t-\Delta w +c_1(x,t)\,w\,=\,u\quad\mbox{a.\,e. in }\,Q,&\\[2mm]
\label{eq:2.9}
&\partial_\nf w + \partial_t w_\Gamma-\dega
w_\Gamma+c_2(x,t)\,w_\Gamma=
\uga, \quad 
\juerg{w_\Gamma=w_{|\Gamma},}
 \quad\mbox{a.\,e. on }\,\Sigma,&\\[2mm]
\label{eq:2.10}
&w(\,\cdot\,,0)=w_0\quad\mbox{a.\,e. in }\,\oma, 
\qquad w_\Gamma(\,\cdot\,,0)=w_{0_\Gamma}\quad\mbox{a.\,e. on }
\,\Gamma,& 
\end{eqnarray}  

{\em has a unique solution $(w,w_\Gamma)\in {\cal Y}$.}

\vspace{1mm}
(ii) \,{\em There exists a constant $\,\widehat{C}>0$, which only depends on $\,\oma$, $T$, $\|c_1\|_\qlio$,
and $\,\|c_2\|_\glisig$, such 
\juerg{that}
 the following holds: whenever $w_0=0$ and $w_{0_\Gamma}=0$ then}
\begin{equation}
\label{eq:2.11}
\|(w,w_\Gamma)\|_{\cal Y}\,\le\,\widehat{C}\,\|(u,u_\Gamma)\|_{\cal H}\,.
\end{equation}
 
\vspace{3mm}
{\em Proof:}  \,\,We put $\beta(w):=\beta_\Gamma(w):=w$ and define the operators 
$$\Pi(w)(x,t):= {(c_1(x,t) - 1) \, w (x,t)},\quad \Pi_\Gamma(w_\Gamma)(x,t):={( c_2(x,t)-1 )\, w_\Gamma(x,t)}.$$
With these definitions, we may rewrite the equations (\ref{eq:2.8}) and (\ref{eq:2.9}) in the form
\begin{eqnarray}
\label{eq:2.12}
&w_t-\Delta w +\beta(w)+\Pi(w)=u,\juerg{\quad\mbox{a.\,e. in }\,Q,}&\\[2mm]
\label{eq:2.13}
&\partial_\nf w+ \partial_t w_\Gamma-\dega w_\Gamma +\beta_\Gamma(w_\Gamma)+\Pi_\Gamma(w_\Gamma)=\uga,\quad \juerg{w_\Gamma=w_{|\Gamma}, \quad\mbox{a.\,e. on }\,\Sigma,} & \qquad
\end{eqnarray}
respectively. Since the functions $\beta$ and $\beta_\Gamma$ are strictly monotone increasing in $\rz$,
the system (\ref{eq:2.10}), (\ref{eq:2.12}), (\ref{eq:2.13}) has almost the same form as the system considered in
Theorem~2.5 in \cite{CC}, the only difference being that \juerg{$\pi, \pi_\Gamma$ in \cite{CC} were Lipschitz continuous functions, while here $\Pi, \Pi_\Gamma$ are linear and (Lipschitz) continuous operators.} However, a closer
inspection of the proof of Theorem~2.5 in \cite{CC} reveals that the argumentation used there carries over
to the present situation with only minor and obvious modifications. \juerg{In particular, we point out that in order to perform the basic a priori estimates leading to the ${\cal Y}$-regularity of the solution (and to estimate (\ref{eq:2.11}) when $w_0=0$ and $w_{0_\Gamma}=0$), equations 
(\ref{eq:2.12})--(\ref{eq:2.13}) (that is, (\ref{eq:2.8})--(\ref{eq:2.9})) are never differentiated in space or time, which is the only source for creating problems in the case of coefficients that depend on $(x,t)$.}
Hence, the asserted existence result is valid.

\vspace{3mm}
\juerg{To show assertion (ii),}  let $\,w_0=0$ and $\,w_{0_\Gamma}=0$. In the following, we denote by $C_i$, $i\in\nz$, positive constants that only depend on the quantities mentioned in the assertion of (ii). {Testing (\ref{eq:2.12}) by} $w_t$ yields for every $t\in (0,T]$ the inequality
\begin{eqnarray*}
&&\txt w_t^2\,\dx\,\dt \,+\,\frac 12\,\|w(t)\|_V^2 \,+\,\txg |\partial_t w_\Gamma|^2\,\dgm\,\dt
\,+\,\frac 12\,\|w_\Gamma(t)\|_{V_\Gamma}^2\\[1mm]
&&\le\,\txt({(|c_1|+1)} |w|\,+\,|u|)\,|w_t|\,\dx\,\dt\,+\,\txg ({(|c_2|+1)}|w_\Gamma|\,+\,|\uga|)\,|\partial_t w_\Gamma|
\,\dgm\,\dt\,,
\end{eqnarray*}
whence, using Young's inequality and the fact that $\,c_1\in\qlio$ and $\,c_2\in\glisig$, we obtain
\begin{eqnarray*}
&&\txt w_t^2\,\dx\,\dt \,+\,\|w(t)\|_V^2 \,+\,\txg |\partial_t w_\Gamma|^2\,\dgm\,\dt
\,+\,\|w_\Gamma{(t)}\|_{V_\Gamma}^2\\[1mm]
&&\le\,C_1\,\Bigl(\txt(|w|^2\,+\,|u|^2)\,\dx\,\dt\,+\,\txg (|w_\Gamma|^2\,+\,|\uga|^2)
\,\dgm\,\dt\Bigr)\,.
\end{eqnarray*}
Gronwall's lemma then yields that
\begin{equation}
\label{eq:2.14}
\|w\|_{H^1(0,T;H)\cap C^0([0,T];V)}\,+\,\|w_\Gamma\|_{H^1(0,T;H_\Gamma)\cap C^0([0,T];V_\Gamma)}\,\le\,
C_2 \,\|(u,u_\Gamma)\|_{\cal H}\,.
\end{equation}

\noindent
Next, a comparison argument in (\ref{eq:2.8}) shows that also
\begin{equation}
\label{eq:2.15}
\|\Delta w\|_{L^2(0,T;H)} \,\le\,C_3 \,\|(u,u_\Gamma)\|_{\cal H}\,.
\end{equation}

\noindent
Now we invoke 
\juerg{\cite[Theorem~3.2, p.~1.79]{Brezzi}}
with the specifications
$$A=-\Delta, \quad g_0=y_{|\Gamma}, \quad p=2, \quad r=0,\quad s=3/2,$$
to conclude that
\begin{equation}
\label{eq:2.16}
\int_0^T\|w(t)\|^2_{H^{3/2}(\oma)}\,\dt \,\le\, C_4\int_0^T\left(\|\Delta w(t)\|_H^2\,+\,
\|w_\Gamma(t)\|_{V_\Gamma}^2\right)\,\dt,
\end{equation}

\noindent
\juerg{so that it follows from (\ref{eq:2.14}) and (\ref{eq:2.15}) 
that} 
\begin{equation}
\label{eq:2.17}
\|w\|_{L^2(0,T;H^{3/2}(\oma))}\,\le\,C_5\,\|(u,u_\Gamma)\|
_{\cal H}\,.
\end{equation}

\noindent
Hence, by the trace theorem \juerg{(cf. \cite[Theorem~2.27, p.~1.64]{Brezzi})}, we have that
\begin{equation}
\label{eq:2.18}
\|\partial_\nf w\|_{L^2(0,T;H_\Gamma)}\,\le\,C_6\,\|(u,u_{\Gamma})\|_{\cal H}\,,
\end{equation}
{so that, by comparison in the equation resulting from (\ref{eq:2.9}), we obtain 
\begin{equation}
\label{eq:2.20}
\|\dega w_\Gamma\|_\glzsig \le\, C_7 \|(u,u_{\Gamma})\|_{\cal H}\,,
\end{equation}
and consequently,
\begin{equation}
\label{eq:2.21}
\|w_\Gamma\|_{L^2(0,T;H^2(\Gamma))} \le\,C_8\,\|(u,u_{\Gamma})\|_{\cal H}\,.
\end{equation}
Finally, owing to standard elliptic estimates, we infer that}
\begin{equation}
\label{eq:2.19}
\|w\|_{L^2(0,T;H^2(\oma))}\,\le\,C_9\,\|(u,u_{\Gamma})\|_{\cal H}\,.
\end{equation}        
This concludes the proof of the assertion.\qed

\vspace{3mm}
{\bf Remark 3:} \,\,It follows from (ii) in Theorem~2.2 that for zero initial data the solution
operator $\uuga\mapsto (w,w_\Gamma)$ is a continuous linear mapping from ${\cal H}$
into ${\cal Y}$.  

\vspace{4mm}
While it cannot be expected that the solution to the linear system (\ref{eq:2.8})--(\ref{eq:2.10})
is bounded, we now establish a boundedness result for the solution to the nonlinear state system
(\ref{eq:1.2})--(\ref{eq:1.4}) that will be 
\juerg{of key importance for}
the subsequent
analysis. To this end, we need the following assumption:

\vspace{2mm}
{\bf (A5)}\quad$y_0\in \lio$, $y_{0_\Gamma} \in L^\infty(\Gamma)$, and it holds
\begin{eqnarray}
\label{eq:2.22}
&&0\,<\,\essinf_{x\in\oma}\,y_0(x), \quad
\esssup_{x\in\oma}\,y_0(x)\,<\,1,\nonumber\\[1mm] 
&&0\,<\,\essinf_{x\in\Gamma}\,y_{0_\Gamma}(x), \quad
\esssup_{x\in\Gamma}\,y_{0_\Gamma}(x)\,<\,1\,. 
\end{eqnarray}

\vspace{2mm}
{\bf Lemma~2.3} \,\,{\em Suppose that the assumptions} {\bf (A2)}--{\bf (A5)} {\em are satisfied. Then
there are constants $0<r_*\le r^*<1$, which only depend on $\Omega$, $T$, $y_0$, $y_{0_\Gamma}$,
$f$, $g$, and $R$, 
such that we have: whenever $(y,y_\Gamma)={\cal S}(u,\uga)$ for some $(u,\uga)\in {\cal U}$ then it holds}
\begin{equation}
\label{eq:2.23}
r_*\le y\le r^* \mbox{\,\, a.\,e. in \,\,}Q, \quad\,r_*\le y_\Gamma\le r^* \mbox{\,\, a.\,e. in \,\,}\Sigma.
\end{equation}

\vspace{2mm}
{\em Proof:} \,\,Let $\uuga\in {\cal U}$ be arbitrary and $(y,y_\Gamma)={\cal S}\uuga$. 
Then we have
$$\|u\|_\qlio\,+\,\|\uga\|_\glisig\,\le\,R.$$
By virtue of (\ref{eq:2.3}) and (\ref{eq:2.22}), there are constants $0<r_*\le r^*<1$ such that
\begin{eqnarray}
\label{eq:2.24}
&&r_*\,\le\,{\rm min}\,\left\{\essinf_{x\in\oma}\,y_0(x)\,,\,\essinf_{x\in\Gamma}\,y_{0_\Gamma}(x)\right\},
\\[1mm]
\label{eq:2.25}
&&r^*\,\ge\,\max\,\left\{\esssup_{x\in\oma}\,y_0(x)\,,\,\esssup_{x\in\Gamma}\,y_{0_\Gamma}(x)\right\},
\\[1mm]
\label{eq:2.26}
&&\max\,\{f'(r)+R\,,\,g'(r)+R\}\,\le\,0 \quad\forall\,r\in (0,r_*),\\[2mm]
\label{eq:2.27}
&&{\rm min}\,\{f'(r)-R\,,\,g'(r)-R\}\,\ge\,0 \quad\forall\,r\in (r^*,1).
\end{eqnarray}
Now define $w:=(y-r^*)^+$. Clearly, we have $w\in V$ and $w_{|\Gamma}\in V_\Gamma$. We put
$w_\Gamma:=w_{|\Gamma}$ and test (\ref{eq:1.2}) by $w$. {Thanks to}
(\ref{eq:2.25}), we readily see that
\begin{eqnarray*}
0&\!=\!&\frac 1 2\,\|w(T)\|_H^2\,+\,\tinto\|\nabla w(t)\|_H^2\,\dt\nonumber\\
&&+\,\frac 1 2 \,\|w_\Gamma(T)\|^2_{H_\Gamma}
\,+\,\tinto\|\nabla_\Gamma w_\Gamma(t)\|_{H_\Gamma}^2\,\dt\,+\,\Phi,
\end{eqnarray*} 
where, owing to (\ref{eq:2.26}) and (\ref{eq:2.27}),
$$
\Phi\,:=\,\txinto (f'(y)-u)\,w\,\dx\,\dt\,+\,\tgamma (g'(y_\Gamma)-\uga)\,w_\Gamma\,\dgm\,\dt
\,\ge\,0\,.
$$
In conclusion, $w=(y-r^*)^+= 0$, i.\,e., $y\le r^*$, almost everywhere in $Q$ and on $\Sigma$. 
The remaining inequalities follow similarly by testing (\ref{eq:1.2}) with $w:=-(y-r_*)^-$.
\qed

\vspace{3mm}
Observe that 
\juerg{on account}
 of {\bf {(A2)}} and of Lemma~2.3, we may (by possibly choosing a larger $K_1^*$)
\juerg{also assume} that 
\begin{equation}
\label{eq:2.28}
\max_{0\le i\le 3}\left\{\max\,\left\{\|f^{(i)}(y)\|_\qlio\,,\,
\|g^{(i)}(\yga)\|_\glisig\right\}\right\}\,\le\,K_1^*\,,
\end{equation}
whenever $(y,\yga)={\cal S}\uuga$ for some $\uuga\in {\cal U}$.
\vspace{2mm}

{\bf Remark 4: }  Lemma~2.3 entails that the singular 
\juerg{parts $\,f_1'\,$ and $\, g_1'\,$ of $\,f'\,$ and $\,g'$, respectively,
which enter the state equations of the control problem 
{\bf (CP)},}
are only active in a domain of arguments where they behave like standard bounded smooth
nonlinearities. As a consequence, we could use classical differentiability results to see that both 
$f$ and $g$ generate three times continuously differentiable Nemytskii operators on suitable
subsets of $\qlio$ and $\glisig$, respectively. From this point, it would in principle be possible to derive the subsequent differentiability results for the control-to-state mapping by using the implicit function theorem. A corresponding approach
was taken in \cite[Chapter~5]{Tr} for the case of standard Neumann boundary conditions not involving dynamic terms or the surface Laplacian. Here, we prefer a direct approach which, while being  
slighthly longer and possibly less elegant than the use of the implicit function theorem, has the advantage of yielding the explicit form of the corresponding derivatives directly.  

\vspace{2mm}
With Lemma~2.3 at hand, we are now able to improve the stability estimate (\ref{eq:2.7}) from Theorem~2.1.

\vspace{2mm}
{\bf Lemma~2.4} \,\,{\em Suppose that} {\bf (A2)}--{\bf (A5)} {\em are 
satisfied. Then there is a constant $K_3^*>0$, which only depends on $\oma$, $T$, $f$, $g$, and $R$,
such that the following holds: whenever $(u_1,u_{1_\Gamma}),  (u_1,u_{2_\Gamma})\in {\cal U}$
are given and $(y_1,y_{1_\Gamma}), (y_2,y_{2_\Gamma})\in {\cal Y}$ are the associated solutions to the state system}
(\ref{eq:1.2})--(\ref{eq:1.4}), {\em then we have }

\vspace{-6mm}
\begin{equation}
\label{eq:2.29}
\|(y_1,y_{1_\Gamma})-(y_2,y_{2_\Gamma})\|_{\cal Y}
\,\le\,K_3^*\,\|(u_1,u_{1_\Gamma})-(u_2,u_{2_\Gamma})\|_{\cal H}\,.
\end{equation}

{\em Proof:} \,\,In the following, we denote by $C_i$, $i\in\nz$, positive constants
that only depend on the quantities mentioned in the assertion. We subtract the state equations (\ref{eq:1.2})--(\ref{eq:1.4}) corresponding to
$((u_i,u_{i_\Gamma}), (y_i,y_{i_\Gamma}))$, $i=1,2$, from each other and multiply the equation resulting
from (\ref{eq:1.2}) by $\partial_t(y_1-y_2)$. Putting $u=u_1-u_2$, $u_\Gamma=u_{1_\Gamma}-u_{2_\Gamma}$, 
$y=y_1-y_2$, and $y_\Gamma=y_{1_\Gamma}-y_{2_\Gamma}$, we have, for
all $t\in [0,T],$
\begin{eqnarray}
\label{eq:2.30}
&&\int_0^t\!\!\xinto y_t^2\,\dx\,\ds\,+\,\frac 1 2\xinto|\nabla
y(t)|^2\,\dx
\,+\,\int_0^t\!\!\int_\Gamma |\partial_t y_\Gamma|^2\,\dgm\,\ds\,+\,\frac 1 2\int_\Gamma
|\nabla_\Gamma y_\Gamma(t)|^2\,\dgm\qquad\nonumber\\
&&\le\,\int_0^t\!\!\xinto |f'(y_1)-f'(y_2)|\,|y_t|\,\dx\,\ds\,+\,\int_0^t\!\!\int_\Gamma
|g'(y_{1_\Gamma})-g'(y_{2_\Gamma})|\,|\partial_t y_\Gamma|\,\dgm\,\ds\nonumber\\
&&\quad + \,\int_0^t\!\!\xinto |u|\,|y_t|\,\dx\,\ds \,+\,\int_0^t\!\!\int_\Gamma
|u_\Gamma|\,|\partial_t y_\Gamma|\,\dgm\,\ds\,.
\end{eqnarray}
Now observe that Lemma~2.3 (see also (\ref{eq:2.28})) and the mean 
value theorem imply that there is some constant $C_1>0$ such that
$$
|f'(y_1)-f'(y_2)|\,\le\,C_1\,|y| \quad\mbox{a.\,e. in }\,Q,\,\quad  |g'(y_{1_\Gamma})-g'(y_{2_\Gamma})|\,\le\,C_1\,
|y_\Gamma|
\quad\mbox{a.\,e. in }\,\Sigma\,.
$$
Hence, it follows from {Young's inequality and (\ref{eq:2.7})} that
\begin{equation}
\label{eq:2.31}
\|y\|_{H^1(0,T;H)\cap C^0([0,T];V)}+\|y_\Gamma\|_{H^1(0,T;H_\Gamma)
\cap C^0([0,T];V_\Gamma)}\,\le\,C_2\,\|(u,u_{\Gamma})\|_{\cal H}\,.
\end{equation}
From this point we may continue as in the proof of {Theorem~2.2} after proving the estimate (\ref{eq:2.14}):
indeed, by the arguments used there, we can repeat the estimates (\ref{eq:2.15}) to (\ref{eq:2.21})
with $(w,w_\Gamma)$ replaced by $(y,\yga)$ to come to the conclusion that (\ref{eq:2.29}) is satisfied.
This concludes the proof of the assertion.\qed

\section{The optimal control problem}
\setcounter{equation}{0}

We now consider the optimal control problem {\bf (CP)} under the general assumptions 
{\bf (A1)}--{\bf {(A4)}}. 
\juerg{Observe that we do not postulate the validity of {\bf (A5)}, 
which means that we cannot
apply Lemma 2.3.}

\subsection{Existence}
We have the following existence result.

\vspace{2mm}
{\bf Theorem~3.1} \,\,{\em Suppose that the general assumptions} {\bf (A1)}--{\bf {(A4)}} {\em are fulfilled.
Then the optimal control problem} {\bf (CP)} {\em admits a solution.}

\vspace{2mm}
{\em Proof:} \,\,Let $\{(u_n, u_{n_\Gamma})\}\subset \uad$ be a minimizing sequence for 
{\bf (CP)}, and
let $(y_n,y_{n_\Gamma})={\cal S}(u_n, u_{n_\Gamma})$, $n\in\nz$. By virtue of the 
\juerg{global estimate (\ref{eq:2.6})},
we may assume (by possibly selecting a suitable subsequence again indexed by $n$) that
there are functions $(\bar u,\bar u_\Gamma)\in {\cal X}$ and $(\bar y,\bar y_\Gamma)\in {\cal Y}$, such that
\begin{eqnarray*}
&u_n\to \bar u& \mbox{weakly-* in }\,\qlio, 
\quad u_{n_\Gamma}\to \bar u_\Gamma \,\,\mbox{weakly-* in }\,\glisig,\\[1mm]
&y_n\to\bar y&\mbox{weakly-* in }\,H^1(0,T;H)\cap L^\infty(0,T;V)\cap L^2(0,T; H^2(\oma))
,\\[1mm]
&y_{n_\Gamma}\to\bar y_\Gamma&\,\,\mbox{weakly-* in }\,H^1(0,T;H_\Gamma)
\cap L^\infty(0,T;V_\Gamma)\cap L^2(0,T;H^2(\Gamma))
\,.
\end{eqnarray*}

\noindent 
In particular, we have
$$\partial_\nf y_n\to\partial_\nf \bar y, 
\quad \dega y_{n_\Gamma}\to\dega\bar y_\Gamma, \quad\,\mbox{both weakly in }\,\glzsig.
$$

\noindent 
Clearly, $(\bar u,\bar u_\Gamma)\in\uad$. Moreover, we infer from standard compact embedding results 
\juerg{(cf.~{\cite[Sect.~8, Cor.~4]{Simon}})}
that, in particular,
\begin{eqnarray}
&y_n\to \bar y&\mbox{strongly in }\,C^0([0,T];H), \label{p1}\\[1mm]
&y_{n_\Gamma}\to \bar y_\Gamma&\,\,\mbox{strongly in }\,C^0([0,T];H_\Gamma). \label{p2}
\end{eqnarray}

\noindent 
But then, by  the Lipschitz continuity of $f'_2$ and $g'_2$ (see {\bf (A2)}),  
also
\begin{eqnarray*}
{f'_2(y_n)\to f'_2(\bar y)}&\mbox{strongly in }\,C^0([0,T];H),\\[1mm]
{g'_2(y_{n_\Gamma})\to g'_2(\bar y_\Gamma)}&\,\,\mbox{strongly in }\,C^0([0,T];H_\Gamma).
\end{eqnarray*}
\juerg{Moreover, (\ref{eq:2.6bis}) and {\bf (A2)} allow us to assume that (possibly
on a subsequence which is again indexed by $n$)}
\begin{eqnarray*}
f'_1(y_n)\to \bar \xi &\mbox{weakly in }\,L^2(0,T;H),\\[1mm]
g'_1(y_{n_\Gamma})\to \bar  \xi_\Gamma&\,\,\mbox{weakly in }\,L^2(0,T;H_\Gamma),
\end{eqnarray*}
for some weak limits $\bar \xi $ and $ \bar \xi_\Gamma.$ 
Since $f_1$ and $g_1$ are convex 
(so that  $f'_1 $  and $ g'_1 $ are increasing),  
the weak convergences above,  along with (\ref{p1})--(\ref{p2}), imply that  $\bar\xi =  f'_1(\bar y)$ 
a.e. in $Q$ and  $\bar\xi_\Gamma = g'_1(\bar y_\Gamma)$ a.e. on $\Sigma$, due to the maximal monotonicity of the operators induced by $f'_1 $ on $L^2(Q)$  and
\juerg{by}
 $ g'_1 $ 
on $L^2(\Sigma)$ (see, e.g., \cite[Prop.~2.5, p.~27]{Brezis}).
At this point,
we may pass to the limit as $n\to\infty$ in the state system (\ref{eq:1.2})--(\ref{eq:1.4})
(written for $(y_n,y_{n_\Gamma}), (u_n,u_{n_\Gamma})$, $n\in\nz$) to conclude that $(\bar y, \bar
y_\Gamma)={\cal S}(\bar u,\bar u_\Gamma)$,
that is, the pair $((\bar u, \bar u_\Gamma), (\bar y,\bar y_\Gamma))$ is admissible for {\bf (CP)}. It then follows
from the \juerg{(weak)} lower sequential semicontinuity of the cost functional $J$ that 
$(\bar u, \bar u_\Gamma)$ is in fact an optimal control for {\bf (CP)}. \qed 

\subsection{Differentiability of the control-to-state operator}
Suppose now that $\buuga\in\uad$ is a local minimizer for {\bf (CP)}, and let $(\bar y,\bar y_\Gamma)={\cal S}
\buuga$ be the associated state. We consider for fixed $(h,h_\Gamma)\in {\cal X}$ the linearized system:

\begin{eqnarray}
\label{eq:3.1}
&\xi_t-\Delta\xi+f''(\bar y)\,\xi=h \quad\mbox{a.\,e. in }\,Q,&\\[1mm]
\label{eq:3.2}
&\partial_\nf \xi+\partial_t\xi_\Gamma-\dega\xi_\Gamma+g''(\bar y_\Gamma)\,\xi_\Gamma=h_\Gamma, \quad\,
\xi_\Gamma=\xi_{|\Gamma}, \,\quad\mbox{a.\,e. on }\,\Sigma,&\\[1mm]
\label{eq:3.3}
&\xi(\,\cdot\,,0)=0 \quad\mbox{a.\,e. in }\,\oma, 
\qquad \xi_\Gamma(\,\cdot\,,0)=0\quad\mbox{a.\,e. on }\,\Gamma.&
\end{eqnarray}
By Theorem~2.2 the system (\ref{eq:3.1})--(\ref{eq:3.3}) admits for every $(h,h_\Gamma)
\in {\cal H}$ (and thus, a fortiori, for every $(h,h_\Gamma)\in{\cal X}$) a unique solution $(\xi,\xi_\Gamma)\in {\cal Y}$, and the linear mapping $(h,h_\Gamma)
\mapsto (\xi,\xi_\Gamma)$ is continuous from ${\cal H}$ into ${\cal Y}$ and thus also from
${\cal X}$ into ${\cal Y}$.

\vspace{2mm}
{\bf Theorem~3.2} \,\,{\em Suppose that the assumptions} {\bf (A2)}--{\bf (A5)} {\em are satisfied. Then we have the following results:}

\vspace{1mm} 
(i) \,\,{\em Let $(\bar u,\bar u_\Gamma)\in {\cal U}$ be arbitrary.  Then the control-to-state mapping 
${\cal S}$, viewed as a mapping from ${\cal X}$ into ${\cal Y}$, is Fr\'echet differentiable at $(\bar u,\bar u_\Gamma)$, and  the  Fr\'echet derivative $D{\cal S}(\bar u,\bar u_\Gamma)$ is given by 
$D{\cal S}(\bar u,\bar u_\Gamma)(h,h_\Gamma) =(\xi,\xi_\Gamma)$, where  for any given
$(h,h_\Gamma)\in {\cal X}$ the pair $(\xi,\xi_\Gamma)$ denotes the solution to
the linearized system} (\ref{eq:3.1})--(\ref{eq:3.3}).

\vspace{1mm}
(ii) \,{\em The mapping $D{\cal S}: {\cal U}\to {\cal L}({\cal X},{\cal Y})$, $\buuga\mapsto D{\cal S}
\buuga$ is Lipschitz continuous on $\cal U$ in the following sense: there is a constant $K_4^*>0$ such that
for all $(u,u_\Gamma),\buuga\in{\cal U}$ and all $(h,h_\Gamma)\in{\cal X}$ it holds}
\begin{equation}
\label{eq:3.4}
\|(D{\cal S}(u,\uga)-D{\cal S}\buuga)(h,h_\Gamma)\|_{\cal Y}\,\le\,K_4^*\,\|(u,\uga)-\buuga\|_{\cal H}
\,\|(h,h_\Gamma)\|_{\cal H}\,.
\end{equation}

\vspace{1mm}
{\em Proof:} \,\,We first show (i). To this end, let $(\bar u,\bar u_\Gamma)\in {\cal U}$ be arbitrarily chosen, and let $(\bar y,\bar y_\Gamma)
={\cal S}(\bar u,\bar u_\Gamma)\in {\cal Y}$ denote the associated solution to the state system. Since ${\cal U}$ is an open
subset of ${\cal X}$, there is some $\lambda>0$ such that for any $(h,h_\Gamma)\in {\cal X}$ with
$\|(h,h_\Gamma)\|_{\cal X}\,\le\,\lambda$ it holds $ (\bar u,\bar u_\Gamma)+(h,h_\Gamma)\in {\cal U}$.
In the following, we consider such variations $(h,h_\Gamma)\in {\cal X}$, and we denote by $(y^h, y_\Gamma^h)$ the
solution to the state system (\ref{eq:1.2})--(\ref{eq:1.4}) associated with $ (\bar u,\bar u_\Gamma)+(h,h_\Gamma)$. 
Moreover, we denote by $(\xi^h,\xi^h_\Gamma)$ the unique solution to the linearized system (\ref{eq:3.1})--(\ref{eq:3.3})
associated with $(h,h_\Gamma)$. We also denote by $C_i$, $i\in\nz$, positive constants that depend neither
on the choice of  $(h,h_\Gamma)\in {\cal X}$ with $\|(h,h_\Gamma)\|_{\cal X}\,\le\,\lambda$ nor on 
$t\in [0,T]$. Now let

\vspace{-6mm}
$$v^h:=y^h-\bar y-\xi^h,\quad\,v^h_\Gamma:=y^h_\Gamma - \bar y_\Gamma-\xi^h_\Gamma.$$
Obviously, we have $(v^h,v^h_\Gamma)\in {\cal Y}$. 
Since the linear mapping $(h,h_\Gamma)\mapsto (\xi^h,\xi_\Gamma^h)$ is by Theorem~2.2~(ii) continuous
from ${\cal X}$ into ${\cal Y}$, it obviously suffices to show that there is an increasing function
$G:[0,\lambda]\to [0,+\infty)$ which satisfies $\lim_{r\searrow 0} \,G(r)/r=0$ and 
\begin{equation}
\label{eq:3.5} 
\|(v^h,v_\Gamma^h)\|_{\cal Y}\,\le\,G\left(\|(h,h_\Gamma)\|_{\cal H}\right)\,.
\end{equation}

\noindent
Apparently, $(v^h, v^h_\Gamma)$ is a solution
to the initial-boundary value problem
\begin{eqnarray}
\label{eq:3.6}
&v^h_t-\Delta v^h\,+\,f'(y^h)-f'(\bar y)-f''(\bar y)\,\xi^h=0 \,\quad\mbox{a.\,e. in }\,Q,&\\[2mm]
\label{eq:3.7}
&\partial_{\bf n}v^h+\partial_t v^h_\Gamma-\Delta_\Gamma
v_\Gamma^h\,+\,g'(y^h_\Gamma)-g'(\bar y_\Gamma)
-g''(\bar y_\Gamma)\,\xi^h_\Gamma=0, \quad \juerg{v^h_\Gamma=v^h_{|\Gamma},}\nonumber\\[1mm]
&\mbox{a.\,e. on }\,\Sigma,&\\[2mm]
\label{eq:3.8}
&v^h(\cdot,0)=0 \quad\mbox{a.\,e. in }\,\oma,\,\,\quad v^h_\Gamma(\cdot,0)=0
\quad\mbox{a.\,e. on }\,\Gamma.&
\end{eqnarray}

Next, we observe that it follows from Taylor's theorem and from the global estimate (\ref{eq:2.23}) {(cf.~also (\ref{eq:2.28}))}
that almost everywhere in $Q$ we have 
\begin{equation}
\label{eq:3.9}
f'(y^h)-f'(\bar y)-f''(\bar y)\,\xi^h\,=\,f''(\bar y)\,v^h\,+\,\Phi^h,
\end{equation}
with some function $\Phi^h\in\qlio$ such that, almost everywhere in $Q$,
\begin{equation}
\label{eq:3.10}
|\Phi^h|\,\le\,\frac 12\,\max_{r_*\le r\le r^*}\,|f^{{(3)}}(r)|\,|y^h-\bar y|^2\,\le\,
\frac {K_1^*} 2\,|y^h-\bar y|^2\,.
\end{equation}

\noindent
By the same token, there is some 
\juerg{$\Phi_\Gamma^h\in\glisig$}
such that, 
almost everywhere on $\Sigma$,

\vspace{-4mm}
\begin{equation}
\label{eq:3.11}
g'(y^h_\Gamma)-g'(\bar y_\Gamma)-g''(\bar y_\Gamma)\,\xi_\Gamma^h\,=\,g''(\bar y_\Gamma)\,v^h_\Gamma\,+\,\Phi_\Gamma^h,
\end{equation}

\vspace{-2mm}
\noindent
where 

\vspace{-6mm}
\begin{equation}
\label{eq:3.12}
|\Phi^h_\Gamma|\,\le\,\frac {K_1^*} 2\,|y^h_\Gamma-\bar y_\Gamma|^2\,.
\end{equation}
Hence, with $c_1:=f''(\bar y)\in \qlio$, $c_2:=g''(\bar y_\Gamma)\in\glisig$, $u:=-\Phi^h\in\qlzo$, and
$u_\Gamma:=-\Phi_\Gamma^h\in L^2(\Sigma)$, we see that the system (\ref{eq:3.6})--(\ref{eq:3.8}) satisfied
by $(v^h,v^h_\Gamma)$ has exactly the structure of the system (\ref{eq:2.8})--(\ref{eq:2.10}). It therefore
follows from (\ref{eq:2.11}) in Theorem~2.2 that
\begin{equation}
\label{eq:3.13} 
\|(v^h,v^h_\Gamma)\|_{\cal Y}\,\le\,C_1\,\|(\Phi^h,\Phi^h_\Gamma)\|_{\cal H}\,.
\end{equation}
Now observe that, owing to the embedding $V\subset L^4(\oma)$ and to
(\ref{eq:2.29}) in Lemma~2.4, 
we have 

\vspace{-6mm}
\begin{eqnarray} 
\label{eq:3.14}
&&\|\Phi^h\|_\qlzo^2\,\le\,C_2\txinto |y^h-\bar y|^4\,\dx\,\dt\,=\,C_2\int_0^T\|y^h(t)-\bar y(t)\|^4_{L^4(\oma)}\,\dt
\nonumber\\[1mm]
&&\le\,
\juerg{C_3}
\,T\,\|y^h-\bar y\|_{C^0([0,T];V)}^4\,\le\,
\juerg{C_4}
\,\|(h,h_\Gamma)\|_{\cal H}^4\,.
\end{eqnarray}
Similar reasoning shows that also
$$
\|\Phi^h_\Gamma\|_\glzsig^2\,\le\,
\juerg{C_5}
\,\|(h,h_\Gamma)\|_{\cal H}^4\,.
$$
In conclusion, (\ref{eq:3.5}) is satisfied for the 
function 
\juerg{$\,G(r)=C_1(\sqrt{C_4}+\sqrt{C_5})\,r^2$},
which concludes the proof of assertion (i).

Next, we show the validity of (ii). To this end, let $\buuga\in {\cal U}$ be arbitrary, and let 
$(k,k_\Gamma)\in{\cal X}$ be such that $(\bar u+k,\bar u_\Gamma+k_\Gamma)\in{\cal U}$. We denote 
$(y^k,y^k_\Gamma)={\cal S}(\bar u+k,\bar u_\Gamma+k_\Gamma)$ and $(\bar y,\bar y_\Gamma)=
{\cal S}\buuga$, and we assume that any $(h,h_\Gamma)\in{\cal X}$ with
$\,\|(h,h_\Gamma)\|_{\cal X}=1\,$ is given. It then suffices to show that there is some $L>0$, independent of $(h,h_\Gamma)$, $\buuga$ and $(k,k_\Gamma)$, such that
\begin{equation}
\label{eq:3.15}
\|(\xi^k,\xi^k_\Gamma)-(\xi,\xi_\Gamma)\|_{\cal Y}\,\le\,L\,\|(k,k_\Gamma)\|_{\cal H}\,,
\end{equation}
where $\,(\xi^k,\xi^k_\Gamma)=D{\cal S}(\bar u+k,\bar u_\Gamma+k_\Gamma)(h,h_\Gamma)\,$ and 
$\,(\xi,\xi_\Gamma)=D{\cal S}\buuga(h,h_\Gamma)$. 
For this purpose, we denote in the following by $K_i$, $i\in\nz$, positive constants that \juerg{depend neither}
on $\buuga$, $(k,k_\Gamma)$, nor on the special choice of $(h,h_\Gamma)\in{\cal X}$ with
$\,\|(h,h_\Gamma)\|_{\cal X}=1\,$.

To begin with, observe that from part (i) it easily follows that $\,(w,w_\Gamma):=
(\xi^k,\xi^k_\Gamma)-(\xi,\xi_\Gamma)\in {\cal Y}\,$ \juerg{satisfies $w_\Gamma=w_{|\Gamma}\,$almost everywhere on $\Sigma$}  and solves the initial-boundary value problem
\begin{eqnarray}
\label{eq:3.16}
\hspace*{-8mm}w_t-\Delta w +f''(\bar y)\,w={{}-\xi^k\,(f''(y^k)-f''(\bar y))} \quad\mbox{a.\,e. in }\,Q,\quad\qquad\\[1mm]
\label{eq:3.17}
\hspace*{-8mm}\partial_\nf w+ \partial_t w_\Gamma-\dega w_\Gamma +g''(\bar y_\Gamma)\,w_\Gamma={{}-\,\xi^k_\Gamma\,(g''(y^k_\Gamma)-g''(\bar y_\Gamma))} \quad\mbox{a.\,e. on }\,\Sigma,\\[1mm]
\label{eq:3.18}
\hspace*{-8mm}w(\,\cdot\,,0)=0 \quad\mbox{a.\,e. in }\,\oma,\qquad w_\Gamma(\,\cdot\,,0)=0\quad\mbox{a.\,e. on }\,
\Gamma.\qquad\qquad\quad
\end{eqnarray}
Hence, it follows from Theorem~2.2 that 
\begin{align}
\label{eq:3.19}
&\|(w,w_\Gamma)\|_{\cal Y}\,\notag\\
&\le\,K_1\,\left(\|\xi^k\,(f''(y^k)-f''(\bar y))\|_\qlzo \,+\,\|\xi^k_\Gamma\,(g''(y^k_\Gamma)-g''(\bar y_\Gamma))\|_\glzsig\right)\,.
\end{align}
Now, by Taylor's theorem and (\ref{eq:2.28}), we have almost everywhere in $Q$ (on $\Sigma$, respectively)
\begin{equation}
\label{eq:3.20}
|f''(y^k)-f''(\bar y)|\,\le \,K_1^*\,|y^k-\bar y|\,\quad\mbox{and }\,\,\, |g''(y^k_\Gamma)-g''(\bar y_\Gamma)|\,\le\,K_1^*\,|y^k_\Gamma-\bar y_\Gamma|\,.
\end{equation}
At this point, we recall that ${\cal U}$ is obviously a bounded subset of ${\cal H}$. Since
$(\bar u+k,\bar u+k_\Gamma)\in{\cal U}$ and $\|(h,h_\Gamma)\|_{\cal X}=1$, we thus can infer from
(\ref{eq:2.28}) and from
the estimate (\ref{eq:2.11}) in Theorem~2.2 that $(\xi^k,\xi^k_\Gamma)$ is bounded in ${\cal Y}$
independently of $(k,k_\Gamma)$, $\buuga$, and the choice of $(h,h_\Gamma)\in{\cal X}$ with
$\|(h,h_\Gamma)\|_{\cal X}=1$. Using the embedding $V\subset L^4(\oma)$ and Lemma~2.4, we 
therefore have
\begin{eqnarray}
\label{eq:3.21}
&&\|\xi^k\,(f''(y^k)-f''(\bar y))\|_\qlzo^2\,\le\,K_2\txinto|\xi^k|^2\,|y^k-\bar y|^2\,\dx\,\dt\nonumber\\[1mm]
&&\le\,K_2\int_0^T\|\xi^k(t)\|^2_{L^4(\oma)}\,\|y^k(t)-\bar y(t)\|^2_{L^4(\oma)}\,\dt\nonumber \\[1mm]
&&\le\,K_3\,\|(y^k-\bar y,y^k_\Gamma-\bar y_\Gamma)\|^2_{\cal Y}
\, \le\,K_4\,\|(k,k_\Gamma)\|_{\cal H}^{{2}}\,.
\end{eqnarray}

Since an analogous estimate holds for the second summand in the bracket on the right-hand side of
(\ref{eq:3.19}), the assertion follows.
\qed

\vspace{3mm}
{\bf Remark 5:} \,\,Notice that we could not establish Fr\'echet differentiability of
${\cal S}\,$ on ${\cal H}$; in fact,
we only were able to show differentiability on the bounded subset ${\cal U}$ of ${\cal X}$. In particular, the 
boundedness of $(u,\uga)$ in ${\cal X}$ was an indispensable prerequisite for proving Lemma~2.3 and the global estimate
(\ref{eq:2.28}), which in turn was fundamental for the derivation of the differentiability result. This will below lead
to a so-called {\em two-norm discrepancy} in the derivation of second-order sufficient optimality
conditions, i.\,e., we will have to work with two different norms.

\vspace{2mm}
With Theorem~3.2 at hand, it is now straightforward to derive the standard variational inequality
that optimal controls must satisfy: indeed, it follows from the quadratic form of $J$ and from the chain rule that the reduced cost functional 
\begin{equation}
\label{eq:3.22}
{\cal J}(u,u_\Gamma):=J({\cal S}(u,u_\Gamma),(u,u_\Gamma))
\end{equation}
is Fr\'echet differentiable at every $\buuga\in {\cal U}$ with the Fr\'echet derivative
\begin{equation}
\label{eq:3.23}
D{\cal J}\buuga=D_{(y,y_\Gamma)}J({\cal S}\buuga,\buuga)\circ D{\cal S}\buuga
\,+\,D_{(u,\uga)}J({\cal S}\buuga,\buuga)\,,
\end{equation}
and, owing to the convexity of $\uad$, we have for every minimizer $\buuga\in\uad$ of ${\cal J}$ 
in $\uad$ that
\begin{equation}
\label{eq:3.24}
D{\cal J}\buuga(v-\bar u,v_\Gamma-\bar u_\Gamma)\,\ge\,0 \quad\forall\,(v,v_\Gamma)\in\uad\,.
\end{equation}
Identification of the expressions in (\ref{eq:3.23}) from (\ref{eq:1.1}) and Theorem~3.2 yields:

\vspace{2mm}
{\bf Corollary 3.3} \,\,{\em Let the assumptions} {\bf (A1)}--{\bf (A5)} {\em be satisfied, and let 
\,$\buuga\in\uad$ be an optimal control for the control problem} {\bf (CP)} {\em with associated 
state $(\bar y,\bar y_\Gamma)={\cal S}\buuga\in {\cal Y}$. 
Then we have, for every $(v,\vga)\in\uad$,}
\begin{eqnarray}    
\label{eq:3.25}
&&\beta_1\txinto (\bar y-z_Q)\,\xi\,\dx\,\dt\,+\,\beta_2\tgamma(\bar y_\Gamma-z_\Sigma)\,\xi_\Gamma\,\dgm\,\dt
\nonumber\\[1mm]
&&+\,\beta_3\xinto (\bar y(\,\cdot\,,T)-z_T)\,\xi(\,\cdot\,,T)\,\dx \,+\,
\beta_4 \int_\Gamma (\bar y_\Gamma(\,\cdot\,,T)-z_{\Gamma,T})\,\xi_\Gamma(\,\cdot\,,T)\,\dgm\nonumber\\[1mm]
&&+\beta_5\int_0^T\!\!
\juerg{\int_\Omega} \bar u(v-\bar u)\,\dx\,\dt
\,+\,\beta_6
\int_0^T\!\!\int_\Gamma \bar u_\Gamma(\vga-\bar u_\Gamma)\,\dgm\,\dt\,\ge\,0\,,            
\end{eqnarray}
{\em where $(\xi,\xi_\Gamma)\in {\cal Y}$ is the unique solution to the linearized system} 
(\ref{eq:3.1})--(\ref{eq:3.3}) {\em associated with $(h,h_\Gamma)=(v-\bar u,\vga-\bar u_\Gamma)$.}

\subsection{First-order necessary optimality conditions}

We are now able to derive the first-order necessary optimality conditions for the control problem {\bf (CP)}.
\juerg{At this point, we have to make the following additional (not overly restrictive)
compatibility assumptions:}

\vspace{3mm}
\juerg{
{\bf (A6)} \quad It holds $\,\beta_3=\beta_4\,$ and $z_{\Gamma,T}=z_{T|\Gamma}$.
}

\vspace{3mm}
\juerg{We then have the following result.}

\vspace{3mm}
{\bf Theorem~3.4} \,\,{\em Let the assumptions} {\bf (A1)}--{\bf 
\juerg {(A6)}} {\em be satisfied, and let 
\,$\buuga\in\uad$ be an optimal control for the control problem} {\bf (CP)} {\em with associated 
state $(\bar y,\bar y_\Gamma)={\cal S}\buuga\in {\cal Y}$. Then the adjoint state system}
\begin{eqnarray}
\label{eq:3.26}
-p_t-\Delta p+f''(\bar y)\,p=\beta_1\,(\bar y-z_Q)\,\quad\mbox{a.\,e. in }\,Q,\qquad\qquad\qquad\\[1mm]
\label{eq:3.27}
\partial_\nf p-\partial_t p_\Gamma-\dega p_\Gamma+g''(\bar y_\Gamma)\,p_\Gamma=\beta_2\,(\bar y_\Gamma-z_\Sigma)\juerg{,\quad p_\Gamma=p_{|\Gamma},}
\,\quad\mbox{a.\,e. on }\,\Sigma,\qquad\\[1mm]
\label{eq:3.28}
p(\,\cdot\,,T)=\beta_3(\bar y(\,\cdot\,,T)-z_T)\quad\mbox{a.\,e. in }\,\oma,
\quad\qquad\qquad\qquad\,\nonumber\\
p_\Gamma(\,\cdot\,,T)=
\juerg{\beta_3}
\,(\bar y_\Gamma(\,\cdot\,,T)-z_{\Gamma,T})\quad\mbox{a.\,e. on }\,\Gamma,\qquad\quad\qquad\,
\end{eqnarray}
{\em has a unique solution $(p,p_\Gamma)\in {\cal Y}$, and for every $(v,v_\Gamma)\in \uad$ we have}

\vspace{-5mm}
\begin{equation}
\label{eq:3.29}
\txinto(p+\beta_5\,\bar u)(v-\bar u)\,\dx\,\dt\,+\,\tgamma (p_\Gamma+\beta_6\,\bar u_\Gamma)(v_\Gamma-
\bar u_\Gamma)\,\dgm\,\dt\,\ge\,0\,.
\end{equation}

\vspace{1mm}
{\em Proof:} \,\,First observe that the system (\ref{eq:3.26})--(\ref{eq:3.28}) is a linear backward-in-time
parabolic initial-boundary value problem, which after the time transformation $t\mapsto T-t$ takes the form
of the system (\ref{eq:2.8})--(\ref{eq:2.10}) provided we put

\vspace{-5mm}
\begin{eqnarray*}
&&c_1(x,t):=f''(\bar y(x,T-t)), \quad c_2(x,t):=g''(\bar y_\Gamma(x,T-t)),\\[1mm]
&&u(x,t):=\beta_1\,(\bar y-z_Q)(x,T-t),\quad
u_\Gamma(x,t):=\beta_2\,(\bar y_\Gamma-z_\Sigma)(x,T-t),\\[1mm]
&&w_0(x):=\beta_3\,(\bar y(x,T)-z_\Gamma(x)),\quad w_{0_\Gamma}(x)=
\juerg{\beta_3}
\,(\bar y_\Gamma(x,T)-z_{\Gamma,T}(x))\,.
\end{eqnarray*}
Obviously, $c_1\in \qlio$, $c_2\in\glisig$, $u\in\qlzo$, $u_\Gamma\in\glzsig$, $w_0\in \heins$, and 
$w_{0_\Gamma}\in\gheins$. 
\juerg{Moreover, it follows from {\bf (A6)} that $\,p(\cdot,T)_{|\Gamma}=
p_\Gamma(\cdot,T)$.}
 Hence, by virtue of Theorem~2.2, the transformed system admits a unique
solution $(w,w_\Gamma)\in {\cal Y}$, so that $(p,p_\Gamma)(x,t):=(w,w_\Gamma)(x,T-t)$ is the unique
solution to the adjoint system, and $(p,p_\Gamma)\in {\cal Y}$.

At this point, we may perform the standard calculation, using repeated integration by parts and
the systems (\ref{eq:3.1})--(\ref{eq:3.3}) and (\ref{eq:3.26})--(\ref{eq:3.28})
\juerg{(for the corresponding calculation in the linear parabolic case with boundary conditions of third kind, we refer the reader to the proof of 
Theorem~3.18 in \cite[p.~158]{Tr}}), which shows that
\begin{eqnarray}
\label{eq:3.30}
&&\beta_1\txinto (\bar y-z_Q)\,\xi\,\dx\,\dt\,+\,\beta_2\tgamma(\bar y_\Gamma-z_\Sigma)\,\xi_\Gamma\,\dgm\,\dt
\nonumber\\
&&+\,\beta_3\xinto (\bar y(\,\cdot\,,T)-z_T)\,\xi(\,\cdot\,,T)\,\dx \,+\,
\juerg{\beta_3}
\int_\Gamma (\bar y_\Gamma(\,\cdot\,,T)-z_{\Gamma,T})\,
\xi_\Gamma(\,\cdot\,,T)\,\dgm\nonumber\\
&&=\txinto p\,h\,\dx\,\dt\,+\,\tgamma p_\Gamma\,h_\Gamma\,\dgm\,\dt\,,
\end{eqnarray}
so  that (\ref{eq:3.29}) follows from  (\ref{eq:3.25}).
\qed

\vspace{3mm}
{\bf Remark 6:} \,\,(i) \,It follows from the above considerations that the Fr\'echet derivative
$D{\cal J}\buuga\in{\cal L}({\cal X},{\cal Y})$ can be identified with the pair $(p+\beta_5\,\bar u,
p_\Gamma+\beta_6\,\bar u_\Gamma)$ in the sense that, with the standard inner product $(\,\cdot\,,\,
\cdot\,)_{\cal H}$ in the Hilbert space ${\cal H}$, we have
\begin{equation}
\label{eq:3.31}
D{\cal J}\buuga(h,h_\Gamma)=((p+\beta_5\,\bar u,
p_\Gamma+\beta_6\,\bar u_\Gamma),(h,h_\Gamma))_{\cal H} \quad\forall\,(h,h_\Gamma)\in {\cal X}\,.
\end{equation}
(ii) \,If $\beta_5>0$ and $\beta_6>0$ then it follows from standard arguments that the condition
(\ref{eq:3.29}) can be given a pointwise interpretation in the following sense:
 we have almost everywhere in $Q$ (on $\Sigma$, respectively) that
$$
\bar u(x,t)=\left\{
\begin{array}{lll}
\widetilde u_2(x,t)&\mbox{if}&\widetilde u_2(x,t)<-\beta^{-1}_5\,p(x,t)\\[1mm]
-\beta^{-1}_5\,p(x,t)&\mbox{if}&\widetilde u_1(x,t)\le -\beta^{-1}_5\,p(x,t)\le \widetilde u_2(x,t)\,,\\[1mm]
\widetilde u_1(x,t)&\mbox{if}&\widetilde u_1(x,t)>-\beta^{-1}_5\,p(x,t)
\end{array}
\right.\qquad\qquad
$$ 

\vspace{-4mm}
\begin{equation}
\bar u_\Gamma(x,t)=\left\{
\begin{array}{lll}
\widetilde u_{2_\Gamma}(x,t)&\mbox{if}&\widetilde u_{2_\Gamma}(x,t)<-\beta^{-1}_6\,p_\Gamma(x,t)\\[1mm]
-\beta^{-1}_6\,p_\Gamma(x,t)&\mbox{if}&\widetilde u_{1_\Gamma}(x,t)\le -\beta^{-1}_6\,p_\Gamma(x,t)
\le \widetilde u_{2_\Gamma}(x,t)\, , \\[1mm]
\widetilde u_{1_\Gamma}(x,t)&\mbox{if}&\widetilde u_{1_\Gamma}(x,t)>-\beta^{-1}_6\,p_\Gamma(x,t)
\end{array}
\right.
\end{equation}
{where $\widetilde u_1, \, \widetilde u_2, \, \widetilde u_{1_\Gamma}, \, 
\widetilde u_{2_\Gamma}$ represent the control constraints defined in {\bf (A1)}.}

\subsection{The second-order Fr\'echet derivative of the control-to-state operator} 

Since the control problem {\bf (CP)} is nonconvex, the first-order necessary optimality conditions established in the previous section are not sufficient. However, it is of utmost importance, e.\,g., for the numerical solution of {\bf (CP)}, to derive sufficient optimality conditions. For this purpose, it is necessary to show that the control-to-state mapping is twice continuously Fr\'echet differentiable. We have the following result. 

\vspace{3mm}
{\bf Theorem~3.5} \,\,\juerg{{\em In addition to} {\bf (A1)}--{\bf (A6)}, {\em assume that}}

\vspace{1mm}
\juerg{{\bf (A7)}}\quad$f,g\in C^4(0,1)$.

\vspace{1mm} \noindent {\em Then we have the following results:}

\vspace{1mm}
(i) \,\,{\em The control-to-state operator ${\cal S}$ is at any $\buuga\in {\cal U}$ twice Fr\'echet
differentiable, and the second Fr\'echet derivative $D^2{\cal S}\buuga\in {\cal L}(\cal X,{\cal L}(\cal X,
\cal Y))$ is defined as follows: if $(h,h_\Gamma), (k,k_\Gamma)\in {\cal X}$ are arbitrary, then
$D^2 {\cal S}\buuga[(h,h_\Gamma)\,,\,(k,k_\Gamma)]=:(\eta,\eta_\Gamma)\in {\cal Y}$ is the unique solution to the 
initial-boundary value problem}   
\begin{eqnarray}
\label{eq:3.33}
\eta_t-\Delta\eta+f''(\bar y)\,\eta=-f^{{(3)}}(\bar y)\,\phi\,\psi \quad\mbox{a.\,e in }\,Q, \qquad\\[1mm]
\label{eq:3.34}
\partial_\nf \eta+\partial_t\eta_\Gamma-\dega\eta_\Gamma+g''(\bar y_\Gamma)\,\eta_\Gamma
=-g^{{(3)}}(\bar y_\Gamma)\,\phi_\Gamma\,\psi_\Gamma
\juerg{, \quad \eta_\Gamma=\eta_{|\Gamma},} \nonumber\\
\quad\mbox{a.\,e. on }\,\Sigma,\\[1mm]
\label{eq:3.35}
\eta(\,\cdot\,,0)=0 \quad\mbox{a.\,e. in }\,\Omega,\qquad \eta_\Gamma(\,\cdot\,,0)=0
\quad\mbox{a.\,e. on }\,\Gamma,\qquad\qquad 
\end{eqnarray}
{\em where we have put}
\begin{equation}
\label{eq:3.36}
(\bar y,\bar y_\Gamma)={\cal S}\buuga, \quad(\phi,\phi_\Gamma)=D{\cal S}(\bar u,\bar  u_\Gamma)(h,h_\Gamma),
\quad (\psi,\psi_\Gamma)=D{\cal S}\buuga(k,k_\Gamma)\,.
\end{equation}

\vspace{2mm}
(ii) \,{\em The mapping $D^2{\cal S}:{\cal U}\to {\cal L}({\cal X},{\cal L}({\cal X},{\cal Y}))$,
$\buuga\mapsto D^2{\cal S}\buuga$, is Lipschitz continuous on ${\cal U}$ in the following sense: there
exists a constant $K_5^*>0$ such that for every $(u,u_\Gamma),\buuga\in{\cal U}$ and all
$(h,h_\Gamma),(k,k_\Gamma)\in{\cal X}$ it holds} 
\begin{eqnarray}
\label{eq:3.37}
&&\|(D^2{\cal S}(u,u_\Gamma)-D^2{\cal S}\buuga)[(h,h_\Gamma),(k,k_\Gamma)]\|_{\cal Y}\nonumber\\[1mm]
&&\le\,K_5^*\,\|(u,u_\Gamma)-(\bar u,\bar u_\Gamma)\|_{\cal H}\,\|(h,h_\Gamma)\|_{\cal H}\,\|(k,k_\Gamma)\|_
{\cal H}\,.
\end{eqnarray}

\vspace{3mm}
{\em Proof:} \,\,We first prove part (i) of the assertion. To this end, let $\buuga\in{\cal U}$ be fixed, $(h,h_\Gamma),
(k,k_\Gamma)\in {\cal X}$ be arbitrary,
and $(\bar y,\bar y_\Gamma),
(\phi,\phi_\Gamma), (\psi,\psi_\Gamma)\in {\cal Y}$ be defined as in (\ref{eq:3.36}). Then, with
\begin{eqnarray*}
&&c_1:=f''(\bar y)\in\qlio, \quad c_2:=g''(\bar y_\Gamma)\in\glisig,\quad u:=-f^{{(3)}}(\bar y)\,\phi\,\psi\in \qlzo,\\[1mm]
&&u_\Gamma:=-g^{{(3)}}(\bar y_\Gamma)\,\phi_\Gamma\,\psi_\Gamma\in\glzsig,
\end{eqnarray*}
the system (\ref{eq:3.33})--(\ref{eq:3.35}) takes the form (\ref{eq:2.8})--(\ref{eq:2.10}) and thus
enjoys a unique solution
$(\eta,\eta_\Gamma)\in {\cal Y}$. Moreover, by (\ref{eq:2.11}) we have the estimate
\begin{equation}
\label{eq:3.38}
\|(\eta,\eta_\Gamma)\|_{\cal Y}\,\le\,\widehat{C}\left( \|f^{{(3)}}(\bar y)\,\phi\,\psi\|_\qlzo\,+\,
\|g^{{(3)}}(\bar y_\Gamma)\,\phi_\Gamma\,\psi_\Gamma\|_\glzsig\right)\,.
\end{equation}

In the remainder of the proof of part (i), we denote by $C_i$, $i\in\nz$, positive constants that do not depend on the
quantities $(h,h_\Gamma)$, $(k,k_\Gamma)$, and $(\bar u,\bar u_\Gamma)$. Using (\ref{eq:2.28}) and 
\juerg{(\ref{eq:3.4})},
and invoking the embedding $V\subset L^4(\Omega)$, we find that
\begin{eqnarray}
\label{eq:3.39}
&&\|f^{{(3)}}(\bar y)\,\phi\,\psi\|_\qlzo^2\,\le\,{K_1^*}^2\txinto |\phi|^2\,|\psi|^2\,\dx\,\dt
\,\le\,C_1\int_0^T \|\phi(t)\|^2_{L^4(\oma)}\,\|\psi(t)\|^2_{L^4(\oma)}\,\dt\nonumber\\[1mm]
&&\le\,C_2\,\|\phi\|^2_{C^0([0,T];V)}\,\|\psi\|^2_{C^0([0,T];V)}\,\le\,C_3\,\|(h,h_\Gamma)\|_{\cal H}^2
\,\|(k,k_\Gamma)\|^2_{\cal H}\,,
\end{eqnarray}
where the validity of the last inequality can be seen as follows: 
by definition (recall (\ref{eq:3.36})),
$(\phi,\phi_\Gamma)$ is the unique solution to the linear problem (\ref{eq:3.1})--(\ref{eq:3.3}) with zero 
initial conditions, and thus we can infer from Theorem 2.2 (see, in particular, (\ref{eq:2.11})) that
$\|(\phi,\phi_\Gamma)\|_{\cal Y}\le\widehat{C}\,\|(h,h_\Gamma)\|_{\cal H}$. By the same token, we conclude
that
$\|(\psi,\psi_\Gamma)\|_{\cal Y}\le\widehat{C}\,\|(k,k_\Gamma)\|_{\cal H}$. The asserted inequality
therefore follows from the definition of the norm of the space 
\juerg{${\cal Y}$. Additionally, we}
 obtain from similar reasoning that also
$$
\|g^{{(3)}}(\bar y_\Gamma)\,\phi_\Gamma\,\psi_\Gamma\|_\glzsig\,\le\,C_4\,\|(h,h_\Gamma)\|_{\cal H}\,
\|(k,k_\Gamma)\|_{\cal H}\,.
$$
In particular, it follows that the bilinear mapping
${\cal X}\times {\cal X}\mapsto {\cal Y}$, $[(k,k_\Gamma),(h,h_\Gamma)]\mapsto(\eta,\eta_\Gamma)$, 
is continuous. 

Now we prove the assertions concerning existence and form of the second Fr\'echet derivative. Since
${\cal U}$ is open, there is some $\lambda>0$ such that $(\bar u+k,\bar u_\Gamma+k_\Gamma)\in {\cal U}$
whenever $\|(k,k_\Gamma)\|_{\cal X}\,\le\,\lambda$. In the following, we only consider such
perturbations $(k,k_\Gamma)\in {\cal X}$. Then for $(y^k,y^k_\Gamma)={\cal S}(\bar u+k,
\bar u_\Gamma+k_\Gamma)$ the global estimates (\ref{eq:2.6}), (\ref{eq:2.28}), (\ref{eq:2.29}) are
satisfied. Without loss of generality, we may also assume that
\begin{equation}
\label{eq:3.40}
\max \,\left\{\|f^{{(4)}}(y^k)\|_\qlio\,,\,\|g^{{(4)}}(y^k_\Gamma)\|_\glisig\right\}\,\le\,K_1^*
\quad\mbox{if }\,\|(k,k_\Gamma)\|_{\cal X}\,\le\,\lambda\,.
\end{equation}
After these preparations, we observe that it suffices to show that
\begin{eqnarray}
\label{eq:3.41}
&&\left\|D{\cal S}(\bar u+k,\bar u_\Gamma+k_\Gamma)-D{\cal S}\buuga-D^2{\cal S}\buuga(k,k_\Gamma)\right\|_{\cal L
({\cal X},{\cal Y})}\nonumber\\[1mm]
&&=\,\sup_{\|(h,h_\Gamma)\|_{\cal X}=1}\,
\left\|\left(D{\cal S}(\bar u+k,\bar u_\Gamma+k_\Gamma)-D{\cal S}\buuga-D^2{\cal S}\buuga(k,k_\Gamma)\right)
(h,h_\Gamma)\right\|_{\cal Y}\nonumber\\[1mm]
&&\le\,G\left(\|(k,k_\Gamma)\|_{\cal H}\right)\,,
\end{eqnarray}
with an increasing function $G:(0,\lambda]\to (0,+\infty)$ that satisfies $\,\lim_{r\searrow 0}\,G(r)/r=0$.

To this end, let $(h,h_\Gamma)\in {\cal X}$ be arbitrary with $\,\|h\|_\qlio\,+\,\|h_\Gamma\|_\glisig=1$.
We put $(\rho,\rho_\Gamma)=D{\cal S}(\bar u+k, 
\juerg{\bar u_\Gamma}+k_\Gamma)(h,h_\Gamma)$, define the pairs $(\phi,\phi_\Gamma),
(\psi,\psi_\Gamma)\in {\cal Y}$ as in (\ref{eq:3.36}), and put 
$$(w,w_\Gamma):=(\rho-\phi-\eta, \rho_\Gamma-\phi_\Gamma-\eta_\Gamma).$$
Then, according to (\ref{eq:3.41}), we need to show that
\begin{equation}
\label{eq:3.42}
\|(w,w_\Gamma)\|_{\cal Y}\,\le\,G\left(\|(k,k_\Gamma)\|_{\cal H}\right)\,.
\end{equation} 
Now, invoking the explicit expressions for the quantities defined above, it is easily seen that 
$(w,w_\Gamma)$ is a solution to the linear
initial-boundary value problem
\begin{eqnarray}
\label{eq:3.43}
w_t-\Delta w+f''(\bar y)\,w=\sigma\quad\mbox{a.\,e. in }\,Q,\quad\qquad
\\[1mm]
\label{eq:3.44}
\partial_\nf w+\partial_t w_\Gamma-\Delta_\Gamma w_\Gamma+g''(\bar y_\Gamma)\,w_\Gamma
=\sigma_\Gamma \juerg{, \quad w_\Gamma=w_{|\Gamma},}\quad\mbox{a.\,e. in }\,\Sigma,
\\[1mm]
\label{eq:3.45}
w(\,\cdot\,,0)=0\quad\mbox{a.\,e. in }\,\oma,\qquad w_\Gamma(\,\cdot\,,0)=0\quad\mbox{a.\,e. on }\,\Gamma,
\end{eqnarray}
where we have put
\begin{eqnarray}
\label{eq:3.46}
&&\sigma:=-\rho\left(f''(y^k)-f''(\bar y)\right)\,+\,f^{{(3)}}(\bar y)\,\phi\,\psi, \nonumber\\[1mm]
&&\sigma_\Gamma:=-\rho_\Gamma\left(g''(y^k_\Gamma)-g''(\bar y_\Gamma)\right)\,+\,g^{{(3)}}(\bar y_\Gamma)
\,\phi_\Gamma\,\psi_\Gamma\,.
\end{eqnarray}
In view of (\ref{eq:2.28}), and since it is easily checked that $(\sigma,\sigma_\Gamma)
\in{\cal H}$, we may again invoke the estimate (\ref{eq:2.11}) in Theorem~2.2 
to conclude that (\ref{eq:3.42}) is satisfied if only 
\begin{equation}
\label{eq:3.47}
\|(\sigma,\sigma_\Gamma)\|_{\cal H}\,\le\,G\left(\|(k,k_\Gamma)\|_{\cal H}\right)\,.
\end{equation} 
Applying Taylor's theorem to $f''$, and recalling (\ref{eq:3.36}), we readily see that there is a function $\,\omega_f\in\qlio$ such that
\begin{equation}
\label{eq:3.48}
f''(y^k)-f''(\bar y)=f^{{(3)}}(\bar y)\,(y^k-\bar y-\psi)\,+\,f^{{(3)}}(\bar y)\,\psi
\,+\,\omega_f\,(y^k-\bar y)^2\quad\mbox{a.\,e. in }\,Q\,.
\end{equation}
Hence, {we have that}
\begin{equation}
\label{eq:3.49}
\sigma\,=-\rho\,{f^{{(3)}}}(\bar y)\,(y^k-\bar y-\psi)\,-\,\psi\,f^{{(3)}}(\bar y)\,(\rho-\phi)\,-\,
\rho\,\omega_f\,(y^k-\bar y)^2\,. 
\end{equation}
Now recall that from the proofs of Fr\'echet differentiablity and of the Lipschitz continuity of the
Fr\'echet derivative (see, in particular, the estimates (\ref{eq:3.4})--(\ref{eq:3.14}) and (\ref{eq:3.15})--(\ref{eq:3.21}), respectively) it follows that
\begin{eqnarray}
\label{eq:3.50}
&&\|(y^k-\bar y-\psi,y^k_\Gamma-\bar y_\Gamma-\psi_\Gamma)\|_{\cal Y}\,\le\,C_1\,\|(k,k_\Gamma)\|_{\cal H}^2\,,
\nonumber\\[1mm]
&&\|(\rho-\phi,\rho_\Gamma-\phi_\Gamma)\|_{\cal Y}\,\le\,C_2\,\|(k,k_\Gamma)\|_{\cal H}\,.
\end{eqnarray}
Moreover, we can infer from Lemma~2.4 that
\begin{equation}
\label{eq:3.51}
\|(y^k-\bar y,y^k_\Gamma-\bar y_\Gamma)\|_{\cal Y}\,\le\,K_3^*\,\|(k,k_\Gamma)\|_{\cal H}\,,
\end{equation}
and it follows from Theorem~2.2 that $\rho$ is bounded in ${\cal Y}$ by a positive 
constant that is independent of $(k,k_\Gamma), (h,h_\Gamma)\in {\cal X}$ with
$\,\|(k,k_\Gamma)\|_{\cal X}\,\le\,\lambda\,$ and $\,\|(h,h_\Gamma)\|_{\cal X}=1$. 
Finally, we conclude from Remark~3 that with a suitable constant $C_3>0$ it holds
\begin{equation}
\label{eq:3.52}
\|(\psi,\psi_\Gamma)\|_{\cal Y}\,\le\,C_3\,\|(k,k_\Gamma)\|_{\cal H}\,.
\end{equation} 
After these preparations, and invoking H\"older's inequality and the embeddings
$V\subset L^4(\oma)$ and $V\subset L^6(\oma)$, we can estimate as follows:
\begin{eqnarray}
\label{eq:3.53}
&&\|\sigma\|_\qlzo^2\,\le\,C_4\txinto\left(|\rho|^2\,|y^k-\bar y-\psi|^2
\,+\,|\psi|^2\,|\rho-\phi|^2\,+\,|\rho|^2\,|y^k-\bar
y|^4\right)\,\dx\,\dt\nonumber\\ 
&&\le\,C_4\!\int_0^T\!\!\!\left(\|\rho(t)\|_{L^4(\oma)}^2\|(y^k-\bar y-\psi)(t)\|_{L^4(\oma)}^2
\,+\,\|\psi(t)\|_{L^4(\oma)}^2\|\rho(t)-\phi(t)\|_{L^4(\oma)}^2\right)\!\dt\nonumber\\
&&\quad +\,C_4\int_0^T\|\rho(t)\|_{L^6(\oma)}^2\,\|y^k(t)-\bar y(t)\|^4_{L^6(\oma)}\,\dt\nonumber\\[1mm]
&&\le\,C_5\,\max_{0\le t\le T}\left(\|\rho(t)\|_V^2\,\|(y^k-\bar y-\psi)(t)\|_V^2\,+\,
\|\psi(t)\|_V^2\,\|\rho(t)-\phi(t)\|_V^2\right.\nonumber\\[1mm]
&&\left. \hspace*{25mm}+\, \|\rho(t)\|_V^2\,\|y^k(t)-\bar y(t)\|_V^4\right)\nonumber\\[1mm]
&&\le\,C_6\,\|(k,k_\Gamma)\|_{\cal H}^4\,.
\end{eqnarray}

By the same reasoning, a similar estimate can be derived for $\|\sigma_\Gamma\|_\glzsig$, which
concludes the proof of the assertion (i).

\vspace{2mm}
Next, we prove the assertion (ii). To this end, suppose that
$\buuga\in{\cal U}$ and {that $(h,h_\Gamma)$ and $(k,k_\Gamma)$  are arbitrarily chosen in $ {\cal X}$}, and let
$(\delta,\delta_\Gamma)\in{\cal X}$ be arbitrary with $(\bar u+\delta,\bar u_\Gamma+\delta_\Gamma)\in
\juerg{{\cal U}}$.
In the following, we will denote by $C_i$, $i\in \nz$, 
positive constants that do not depend on any of these quantities.
We put
\begin{eqnarray*}
&&(y^\delta,y^\delta_\Gamma)={\cal S}(\bar u+\delta,\bar u_\Gamma+\delta_\Gamma),\quad (\bar y,\bar y_\Gamma)={\cal S}\buuga,
\quad (\phi,\phi_\Gamma)=D{\cal S}\buuga(h,h_\Gamma),\\[1mm]
&&(\psi,\psi_\Gamma)=D{\cal S}\buuga(k,k_\Gamma),\quad (\phi^\delta,
\phi^\delta_\Gamma)=D{\cal S}(\bar u+\delta,\bar u_\Gamma+\delta_\Gamma)(h,h_\Gamma),\qquad\qquad\qquad\\[1mm]
&&(\psi^\delta,\psi^\delta_\Gamma)= D{\cal S}(\bar u+\delta,\bar u_\Gamma+\delta_\Gamma)(k,k_\Gamma) ,\quad
(\eta,\eta_\Gamma)=D^2{\cal S}\buuga[(h,h_\Gamma),(k,k_\Gamma)],\nonumber\\[1mm]
&&(\eta^\delta,\eta_\Gamma^\delta)= D^2{\cal S}(\bar u+\delta,\bar u_\Gamma+\delta_\Gamma)
[(h,h_\Gamma),(k,k_\Gamma)]\,.
\end{eqnarray*}

From the previous results, in particular, (\ref{eq:2.29}) and (\ref{eq:3.4}),
we can infer that there is a constant $C_1>0$ such that
\begin{eqnarray}
\label{eq:3.54}
&&\|(\phi,\phi_\Gamma)\|_{\cal Y}+\|(\phi^\delta,\phi^\delta_\Gamma)\|_{\cal Y}\,\,\le\,C_1\,\|(h,h_\Gamma)\|_{\cal H},
\nonumber\\[1mm]
&&\|(\psi,\psi_\Gamma)\|_{\cal Y}+\|(\psi^\delta,\psi^\delta_\Gamma)\|_{\cal Y}\,\le\,C_1\,\|(k,k_\Gamma)\|_{\cal H},\nonumber\\[1mm]
&&\|(\eta,\eta_\Gamma)\|_{\cal Y}\,\,+\,\|(\eta^\delta,\eta^\delta_\Gamma)\|_{\cal Y}\,\,\le\,C_1\,\|(h,h_\Gamma)\|_{\cal H}\,\|(k,k_\Gamma)\|_{\cal H},\nonumber\\[1mm]
&&\|(y^\delta,y^\delta_\Gamma)-(\bar y,\bar y_\Gamma)\|_{\cal Y}\,\le\,C_1\,\|(\delta,\delta_\Gamma)\|_{\cal H},
\nonumber\\[1mm]
&&\|(\phi^\delta,\phi^\delta_\Gamma)-(\phi,\phi_\Gamma)\|_{\cal Y}\,\le\,C_1\,\|(\delta,\delta_\Gamma)\|_{\cal H}\, 
\|(h,h_\Gamma)\|_{\cal H},\nonumber\\[1mm]
&&\|(\psi^\delta,\psi^\delta_\Gamma)-(\psi,\psi_\Gamma)\|_{\cal Y}\,\le\,C_1\,
\|(\delta,\delta_\Gamma)\|_{\cal H}\,\|(k,k_\Gamma)\|_{\cal H}\,.
\end{eqnarray}

Now observe that $(w,w_\Gamma)=(\eta^\delta,\eta^\delta_\Gamma)-(\eta,\eta_\Gamma)$ satisfies the linear
initial-boundary value problem of the type (\ref{eq:2.8})--(\ref{eq:2.10})
\begin{eqnarray}
\label{eq:3.55}
w_t-\Delta w+f''(\bar y)\,w=\sigma\quad\mbox{a.\,e. in }\,Q,\qquad\qquad\\[1mm]
\label{eq:3.56}
\partial_\nf w+\partial_t w_\Gamma-\dega w_\Gamma+g''(\bar y_\Gamma)\,w_\Gamma=\sigma_\Gamma
\juerg{, \quad w_\Gamma=w_{|\Gamma},}
\quad\mbox{a.\,e. on }\,\Sigma,\,\,\\[1mm]
\label{eq:3.57}
w(\,\cdot\,,0)=0 \quad\mbox{a.\,e. in }\,\oma,\qquad w_\Gamma(\,\cdot\,,0)=0\quad\mbox{a.\,e. on }\,\Gamma,
\end{eqnarray}
where we have put
\begin{eqnarray}
\label{eq:3.58}
&&\sigma=-\eta^\delta(f''(y^\delta)-f''(\bar y))-(f^{{(3)}}(y^\delta)\,\phi^\delta\,\psi^\delta
-f^{{(3)}}(\bar y)\,\phi\,\psi)\nonumber\,,\\[1mm]
&&\sigma_\Gamma=-\eta^\delta_\Gamma
(g''(y^\delta_\Gamma)-g''(\bar y_\Gamma))
-(g^{{(3)}}(y_\Gamma^\delta)\,\phi_\Gamma^\delta\,\psi_\Gamma^\delta-g^{{(3)}}(\bar y_\Gamma)\,\phi_\Gamma\,\psi_\Gamma)\,.
\end{eqnarray}
From Theorem~2.2 it follows that
\begin{equation}
\label{eq:3.59}
\|(w,w_\Gamma)\|_{\cal Y}\,\le\,\widehat C\,\|(\sigma,\sigma_\Gamma)\|_{\cal H}\,,
\end{equation}
so that it remains to show an estimate of the form
\begin{eqnarray}
\label{eq:3.60}
&&\|(\sigma,\sigma_\Gamma)\|_{\cal H}\,\le\,C_2\,\|(\delta,\delta_\Gamma)\|_{\cal H}\,\|(h,h_\Gamma)\|_{\cal H}\,\|(k,k_\Gamma)\|_{\cal H}\,.
\end{eqnarray}
Moreover, we can infer from 
\juerg{(\ref{eq:2.28}), (\ref{eq:2.29}), (\ref{eq:3.4}), and (\ref{eq:3.40}),}
 that, almost everywhere in $Q$, 
\begin{equation}
\label{eq:3.61}
|\sigma|\,\le\,K_1^*\,(|\eta^\delta|\,|y^\delta-\bar y|\,+\,|\phi^\delta|\,|\psi^\delta|\,
|y^\delta-\bar y|\,+\,|\phi^\delta|\,|\psi^\delta-\psi|\,+\,|\psi|\,|\phi^\delta-\phi|)\,.
\end{equation}
Hence, by (\ref{eq:3.54}), and using H\"older's inequality and the embedding $V\subset L^4(\oma)$,{
\begin{align} \label{eq:3.62}
\txinto|\eta^\delta|^2\,|y^\delta-\bar y|^2\,\dx\,\dt \,
&\le\,\int_0^T \|\eta^\delta(t)\|^2_{L^4(\oma)} \,
\juerg{\|(y^\delta-\bar y)(t)\|_{L^4(\oma)}^2}
\,\dt\nonumber\\[1mm]
&\le\,C_3\,\|\eta^\delta\|_{C^0([0,T];V)}^2\,\|y^\delta-\bar y\|^2_{C^0([0,T];V)}
\nonumber\\[2mm]
&\le\,
C_4\,\|(\delta,\delta_\Gamma)\|_{\cal H}^2\,\|(h,h_\Gamma)\|_{\cal H}^2
\,\|(k,k_\Gamma)\|^2_{\cal H}\,.
\end{align}
}%
Similar reasoning yields
\begin{equation}
\label{eq:3.63}
\|\phi^\delta(\psi^\delta-\psi)\|_\qlzo^2\,+\,\|\psi(\phi^\delta-\phi)\|_\qlzo^2\,
\le\,C_5\,\|(\delta,\delta_\Gamma)\|_{\cal H}^2\,\|(h,h_\Gamma)\|_{\cal H}^2
\,\|(k,k_\Gamma)\|_{\cal H}^2\,.
\end{equation}
Moreover, we invoke (\ref{eq:3.54}), H\"older's inequality, and the embeddings $V\subset L^4(\oma)$
and $H^2(\oma)\subset\lio$, to conclude that
\begin{align}
\label{eq:3.64}
&\txinto \!\!|\phi^\delta|^2\,|\psi^\delta|^2\,|y^\delta-\bar y|^2\,\dx\,\dt\le
\int_0^T\!\!\!\|(y^\delta-\bar y)(t)\|_\lio^2\,
\juerg{\|\phi^\delta(t)\|^2_{L^4(\oma)}}
\,\|\psi^\delta(t)\|^2_{L^4(\oma)}\,\dt\nonumber\\[1mm]
&\le\,C_6\,\|\phi^\delta\|^2_{C^0([0,T];V)}\,\|\psi^\delta\|^2_{C^0([0,T];V)}\,\|y^\delta-\bar y\|^2_{L^2(0,T;H^2(\oma))}\hspace*{4cm}\nonumber\\[2mm]
&\le\,C_7\,\|(\delta,\delta_\Gamma)\|_{\cal H}^2\,\|(h,h_\Gamma)\|_{\cal H}^2
\,\|(k,k_\Gamma)\|_{\cal H}^2\,.
\end{align}
Finally, we can derive estimates for $\|\sigma_\Gamma\|_\glzsig$  similar to 
(\ref{eq:3.61})--(\ref{eq:3.64}), which proves the validity of the 
required estimate (\ref{eq:3.60}).
The assertion is thus  proved.
\qed

\subsection{Second-order sufficient optimality conditions}
\juerg{In the following, we assume that the assumptions {\bf (A1)}--{\bf (A7)}
are satisfied.}
With Theorem~3.5 at hand, the road is paved to derive sufficient conditions for
optimality. But, because the control-to-state operator ${\cal S}$ is not Fr\'echet differentiable on
${\cal H}$, we are faced with the two-norm discrepancy, which makes it impossible to establish second-order sufficient optimality conditions by means of the same simple arguments as in the finite-dimensional case or, e.\,g.,
in the proof of Theorem~4.23 on page 231 in \cite{Tr}. It will thus be necessary to tailor the conditions in such a way as to
overcome the two-norm discrepancy. At the same time, for practical purposes the conditions should not be 
overly restrictive. For such an approach, we follow the lines of Chapter 5 in \cite{Tr}, here. Since many of the 
arguments developed here are rather similar to those employed in \cite{Tr}, we can afford to be sketchy 
and refer the reader to \cite{Tr} for
full details.

To begin with, the quadratic cost functional $J$, viewed as a mapping on ${\cal Y}\times{\cal U}$, is obviously
twice continuously Fr\'echet differentiable on ${\cal Y}\times{\cal U}$, and for any 
$((\bar y,\bar y_\Gamma),\buuga)\in {\cal Y}\times{\cal U}$ and any $((v,v_\Gamma),(h,h_\Gamma)),
((w,w_\Gamma),(k,k_\Gamma))\in {\cal Y}\times {\cal X}$ it holds

\vspace{-6mm}
\begin{eqnarray}
\label{eq:3.65}
&&D^2J((\bar y,\bar y_\Gamma),\buuga)[((v,v_\Gamma),(h,h_\Gamma)),
((w,w_\Gamma),(k,k_\Gamma))]\nonumber\\[1mm]
&&=\beta_1\txinto v\,w\,\dx\,\dt\,+\,\beta_2\tgamma v_\Gamma\,w_\Gamma\,\dgm\,\dt\,+\,\beta_3\xinto v(T)\,w(T)\,\dx\nonumber\\
&&\quad+\,
\juerg{\beta_3}
\int_\Gamma v_\Gamma(T)\,w_\Gamma(T)\,\dgm\,+\,\beta_5\txinto h\,k\,\dx\,\dt
\,+\,\beta_6\tgamma h_\Gamma\,k_\Gamma\,\dgm\,\dt\,.\qquad
\end{eqnarray}
Hence, it follows from Theorem~3.5 and from the chain rule that the reduced cost functional ${\cal J}$
is twice continuously Fr\'echet differentiable on ${\cal U}$. Now let $\buuga\in{\cal U}$ be fixed
and $(h,h_\Gamma),(k,k_\Gamma)\in{\cal X}$ be arbitrary. As previously,
we put

\vspace{-5mm}
\begin{eqnarray*}
&&(\bar y,\bar y_\Gamma)={\cal S}\buuga,\quad (\phi,\phi_\Gamma)=D{\cal S}\buuga(h,h_\Gamma),
\quad (\psi,\psi_\Gamma)=D{\cal S}\buuga(k,k_\Gamma),\\[1mm]
&&(\eta,\eta_\Gamma)=D^2{\cal S}\buuga [(h,h_\Gamma),(k,k_\Gamma)]\,.
\end{eqnarray*}
Then a straightforward calculation resembling that carried out on page 241 in \cite{Tr}, using the 
chain rule as main tool, yields the equality

\vspace{-5mm}
\begin{eqnarray}
\label{eq:3.66}
&&D^2{\cal J}\buuga[(h,h_\Gamma),(k,k_\Gamma)]\,=\,D_{(y,y_\Gamma)}J((\bar y,\bar y_\Gamma),\buuga)(\eta,\eta_\Gamma)\nonumber\\[1mm]
&&+\,D^2J((\bar y,\bar y_\Gamma),\buuga)
[((\phi,\phi_\Gamma),(h,h_\Gamma))\,,\,((\psi,\psi_\Gamma),(k,k_\Gamma))]\,.
\end{eqnarray} 
Now observe that the first summand of the right-hand side of (\ref{eq:3.66}) is equal to the expression

\vspace{-6mm}
\begin{eqnarray}
\label{eq:3.67}
&&\beta_1\txinto (\bar y-z_Q)\,\eta\,\dx\,\dt\,+\,\beta_2\tgamma(\bar y_\Gamma-z_\Sigma)\,\eta_\Gamma\,\dgm\,\dt\nonumber\\
&&+\,\beta_3\xinto(\bar y(T)-z_T)\,\eta(T)\,\dx\,+\,
\juerg{\beta_3}
\int_\Gamma(y_\Gamma(T)-z_{\Gamma,T})\,
\eta_\Gamma(T)\,\dgm 
\end{eqnarray} 
and that $(\eta,\eta_\Gamma)$ solves a system of the 
form (\ref{eq:3.1})--(\ref{eq:3.3}), where $h$ is replaced
by $\,-f^{{(3)}}(\bar y)\phi\psi\in \qlzo$ and $h_\Gamma$ by $\,-g^{{(3)}}(\bar y_\Gamma)\phi_\Gamma\psi_\Gamma
\in \glzsig$. Since the calculation leading to (\ref{eq:3.30}) also works for right-hand sides
in $\qlzo\times\glzsig$, we can infer that  

\vspace{-5mm}
\begin{eqnarray}
\label{eq:3.68}
&&D_{(y,y_\Gamma)}J((\bar y,\bar y_\Gamma),\buuga)(\eta,\eta_\Gamma)\nonumber\\
&&=\,-\txinto p\,f^{{(3)}}(\bar y)\,\phi\,\psi\,\dx\,\dt\,-\,\tgamma p_\Gamma\,g^{{(3)}}(\bar y_\Gamma)\,\phi_\Gamma
\,\psi_\Gamma\,\dgm\,\dt\,,
\end{eqnarray}
where $(p,p_\Gamma)\in {\cal Y}$ is the adjoint state associated with $((\bar y,\bar y_\Gamma),\buuga)$.
Summarizing, we have thus shown that it holds the representation formula
\begin{eqnarray}
\label{eq:3.69}  
&&D^2{\cal J}\buuga[(h,h_\Gamma),(h,h_\Gamma)]\,=\,\txinto(\beta_1-p\,f^{{(3)}}(\bar y))\,|\phi|^2\,\dx\,\dt
\nonumber\\
&&+\tgamma(\beta_2-p_\Gamma\,g^{{(3)}}(\bar y_\Gamma))\,|\phi_\Gamma|^2\,\dgm\,\dt
\,+\,\beta_3\xinto|\phi(T)|^2\,\dx \,+\,
\juerg{\beta_3}
\int_\Gamma|\phi_\Gamma(T)|^2\,\dgm\nonumber\\[1mm]
&&+\,\beta_5\,\|h\|^2_\qlzo\,+\,\beta_6\,\|h_\Gamma\|_\glzsig^2\,.
\end{eqnarray}

Equality (\ref{eq:3.69}) gives rise to hope that, under appropriate conditions, $\,D^2{\cal J}\buuga\,$ 
might be a positive definite operator on a suitable subset of the space ${\cal H}$.  
To formulate such a condition, assume that $\buuga\in\uad$ is a given control with associated state
$(\bar y,\bar y_\Gamma)={\cal S}\buuga\in{\cal Y}$ and adjoint state 
\juerg{$(p,p_\Gamma)\in {\cal Y}$}
 satisfying 
(\ref{eq:3.26})--(\ref{eq:3.28}). We then introduce
for fixed $\tau>0$ the {\em set of strongly active constraints for $\buuga$} by
\begin{eqnarray}
\label{eq:3.70}
A_\tau\buuga&\!\!:=\!\!&\{(x,t)\in Q:\,|p(x,t)+\beta_5\,\bar u(x,t)|>\tau\}\nonumber\\
&&\cup \,\{(x,t)\in \Sigma:\,|p_\Gamma(x,t)+\beta_6\,\bar u_\Gamma(x,t)|>\tau\}\,.
\end{eqnarray}   
Apparently it follows from (\ref{eq:3.29}) that, depending on the signs of $p(x,t)+\beta_5\,\bar u(x,t)$ and
of $p_\Gamma(x,t)+\beta_6\,\bar u_\Gamma(x,t)$, the control values $\bar u(x,t)$ and $\bar u_\Gamma(x,t)$,
respectively, attain one of the constraint values. We now define the 
{\em $\tau-$critical cone} $C_\tau\buuga$ to be the set of all $(h,
h_\Gamma)\in{\cal X}$
such that
\begin{eqnarray}
\label{eq:3.71}
&&\,\,h(x,t)\left\{
\begin{array}{lll}
=0&\mbox{if}&(x,t)\in A_\tau\buuga\\[1mm]
\ge 0&\mbox{if}&\bar u(x,t)=\widetilde u_1(x,t) \mbox{\, and }\,(x,t)\not \in A_\tau\buuga\\[1mm]
\le 0&\mbox{if}&\bar u(x,t)=\widetilde u_2(x,t) \mbox{\, and }\,(x,t)\not \in A_\tau\buuga
\end{array}
\right. \,,\nonumber\\[3mm]
&&h_\Gamma(x,t)\left\{
\begin{array}{lll}
=0&\mbox{if}&(x,t)\in A_\tau\buuga\\[1mm]
\ge 0&\mbox{if}&\bar u_\Gamma(x,t)=\widetilde u_{1_\Gamma}(x,t) \mbox{\, and }\,(x,t)\not \in A_\tau\buuga\\[1mm]
\le 0&\mbox{if}&\bar u_\Gamma(x,t)=\widetilde u_{2_\Gamma}(x,t) \mbox{\, and }\,(x,t)\not \in A_\tau\buuga
\end{array}
\right. \,.\qquad
\end{eqnarray}

After these preparations, we can formulate 
the second-order sufficient optimality condition as follows:
\begin{eqnarray}
\label{eq:3.72}
&&\mbox{There exist constants }\,\delta>0\,\mbox{ and
}\,\tau>0\,\mbox{ such that }\nonumber\\
&&\qquad\qquad D^2{\cal J}[(\bar u,\bar u_\Gamma)\left[(h,h_\Gamma),(h,h_\Gamma)\right]\,\ge\,\delta\,\|(h,h_\Gamma)\|_{\cal H}^2\, \quad\forall\,(h,h_\Gamma)\in C_\tau\buuga\,,\qquad\nonumber\\[1mm]
&&\mbox{where }\,D^2{\cal J}[(\bar u,\bar u_\Gamma)\left[(h,h_\Gamma),(h,h_\Gamma)\right]\,
\mbox{ is given by (\ref{eq:3.69})
with }\,(\bar y,\bar y_\Gamma)={\cal S}\buuga,\nonumber\\
&&(\phi,\phi_\Gamma)=D{\cal S}\buuga(h,h_\Gamma) \,\mbox{ and the associated adjoint state }\,(p,p_\Gamma)\,. 
\end{eqnarray}

The following result resembles Theorem~5.17 in \cite{Tr}.

\vspace{2mm}
{\bf Theorem~3.6} \,\,{\em Suppose that the assumptions} {\bf (A1)}--{\bf 
\juerg{(A7)}}
 {\em are satisfied, and
assume that the triple $\buuga\in\uad$, $(\bar y ,\bar y_\Gamma)={\cal S}\buuga\in{\cal Y}$ and 
$(p,p_\Gamma)\in{\cal Y}$ fulfills the first-order necessary optimality conditions}
(\ref{eq:3.26})--(\ref{eq:3.29}). {\em Moreover, assume that the condition} (\ref{eq:3.72}) {\em 
is fulfilled. Then there are constants $\varepsilon>0$ and $\sigma>0$ such that
\begin{eqnarray}
&&{\cal J}(u,u_\Gamma)\,\ge\,{\cal J}\buuga\,+\,\sigma\,\|(u-\bar u,\uga-\bar u_\Gamma)\|^2_{\cal H}\nonumber\\
&&\mbox{whenever }\,(u,\uga)\in\uad\,\mbox{ and }\,\|(u,u_\Gamma)-\buuga\|_{\cal X}\,\le\,\varepsilon\,.
\end{eqnarray}
In particular, $\buuga$ is locally optimal in the sense of ${\cal X}$.}

\vspace{1mm}
{\em Proof:} \,\,The proof closely resembles that of \cite[Theorem~5.17]{Tr}, and therefore we can refer to
\cite{Tr}. We only indicate one argument that needs a bit more explanation. To this end, let  
$(u,\uga)\in\uad$ be arbitrary.  Since ${\cal J}$ is twice continuously Fr\'echet 
differentiable in ${\cal U}$,
it follows from Taylor's theorem with integral remainder (see, e.\,g., 
Theorem~8.14.3 on page 186 in \cite{D}) that  
\begin{eqnarray}
\label{eq:3.73}
{\cal J}(u,\uga)-{\cal J}\buuga&=&DJ\buuga(v,\vga)\,+\,\frac 12\, D^2{\cal J}\buuga[(v,\vga),(v,\vga)]
\nonumber\\
&&+\,R^{\cal J}((u,\uga),\buuga)\,,
\end{eqnarray}
with \juerg{$\,(v,v_\Gamma):=(u-\bar u,\uga-\bar u_\Gamma)\,$ and the remainder} 
\begin{eqnarray}
\label{eq:3.74}
\hspace*{-8mm}&&R^{\cal J}((u,\uga),\buuga)\nonumber\\
\hspace*{-8mm}&&=\,\int_0^1(1-s)\left(D^2{\cal J}(\bar u+sv,\bar u_\Gamma+sv_\Gamma)
-D^2{\cal J}\buuga\right)\left[(v,\vga),(v,\vga)\right]\ds\,.
\end{eqnarray}
A lengthy but straightforward calculation, based on the representation formulas (\ref{eq:3.65})--(\ref{eq:3.67}) as well as on the Lipschitz estimates (\ref{eq:2.29}), (\ref{eq:3.4}), and
(\ref{eq:3.37}), reveals that 
\begin{eqnarray}
\label{eq:3.75}
\left|R^{\cal J}((u,\uga),\buuga)\right|&\le&C_1\int_0^1(1-s)\,s\,\|(v,v_\Gamma)\|_{\cal H}^3\ds\nonumber\\
&\le &C_2\,\|(v,v_\Gamma)\|_{\cal X}\,\|(v,v_\Gamma)\|_{\cal H}^2\,,
\end{eqnarray}
with constants $C_1>0$ and $C_2>0$ that do not depend on the choice of $(u,\uga)\in\uad$.
 From this point, we can argue along exactly the same lines as on pages~292--294 in the proof of
\cite[Theorem~5.17]{Tr} to conclude the validity of the assertion.  
\qed

\vspace{7mm}
\section*{Acknowledgements}
\juerg{The authors gratefully acknowledge fruitful 
discussions with M.~Hassan Farshbaf-Shaker (Weierstrass Institute, Berlin) and the valuable comments of an anonymous referee that helped to improve the readability of the present paper. The warm hospitality of the IMATI of CNR in Pavia 
and of the WIAS in Berlin has also been  appreciated by JS and PC, respectively. Some financial support comes from the FP7-IDEAS-ERC-StG Grant \#200947 (BioSMA) 
and from the MIUR-PRIN Grant 2010A2TFX2 ``Calculus of Variations''.}

\end{document}